\documentclass{amsart}
\usepackage[T1]{fontenc}
\usepackage{latexsym,amssymb,amsmath,mathrsfs}

\usepackage{amsthm}
\usepackage{longtable}
\usepackage{epsfig}
\usepackage{hhline}
\usepackage{epic}

   \newcommand{\cyc}[1]{\langle\,#1\,\rangle}
   
   \newcommand{\Irr}{\operatorname{Irr}}

   \newcommand{\N}{\mathbb{N}}
   
   \newcommand{\Z}{\mathbb{Z}}

   \newcommand{\la}{\lambda}
      \newcommand{\tla}{\tilde{\lambda}}
 \newcommand{\modp}{ \; (\mbox{mod} \; p)}
  \newcommand{\modquat}{ \; (\mbox{mod} \; 4)}
   \newcommand{\modtwo}{ \; (\mbox{mod} \; 2)}
 \newcommand{\da}{\delta}
  \newcommand{\ba}{\beta}
  \newcommand{\sa}{\sigma}
    
     \newcommand{\za}{\zeta}
      \newcommand{\ta}{\theta}
 \newcommand{\ga}{\gamma}
 \newcommand{\e}{\varepsilon}
 
  \newcommand{\p}{\bar{p}}
  \newcommand{\h}{\bar{h}}

\newcommand{\cal}[1]{\mathcal{#1}}

\newcommand{\cla}{\cyc{\la}}
\newcommand{\ctla}{\cyc{\tla}}
 \newcommand{\cQi}{\cyc{Q_i}}

  \newcommand{\floor}[1]{\lfloor #1 \rfloor}

\newtheorem{theorem}{Theorem}[section] % 1st argument is your name for it
\newtheorem{lemma}[theorem]{Lemma}     % 2nd argument is what is printed

\newtheorem{conjecture}[theorem]{Conjecture}
\newtheorem{proposition}[theorem]{Proposition}

\title[Isaacs-Navarro Conjecture for covering groups]% end with percent
   {The Isaacs-Navarro Conjecture for covering groups of the symmetric and alternating groups in odd characteristic}

\author{Jean-Baptiste Gramain}

\address{Department of Mathematical sciences\\
University of Copenhagen\\
Universitetsparken 5\\
DK-2100 Copenhagen O\\}
\email{gramain@math.ku.dk}

\begin{document}

\begin{abstract}
In this paper, we prove that a refinement of the Alperin-McKay Conjecture for $p$-blocks of finite groups, formulated by I. M. Isaacs and G. Navarro in 2002, holds for all covering groups of the symmetric and alternating groups, whenever $p$ is an odd prime.
\end{abstract}

\maketitle

\section{Introduction}\label{intro}

I.M. Isaacs and G. Navarro have formulated in \cite{Isaacs-Navarro} some refinements of the McKay and Alperin-McKay Conjectures for arbitrary finite groups. Consider a finite group $G$ and a prime $p$. Let $B$ be a $p$-block of $G$, with defect group $D$, and let $b$ be the Brauer correspondent of $B$ in $N_G(D)$. Throughout this paper, we will use a $p$-valuation $\nu$ on $\Z$, given by $\nu(n)=a$ if $n=p^a q$ with $(p, \, q)=1$. The height ${\mathfrak{h}}(\chi) \in \Z_{\geq 0}$ of an irreducible (complex) character $\chi \in B$ is then defined by the equality $\nu (\chi(1))=\nu(|G|)-\nu(|D|)+{\mathfrak{h}}(\chi)$. We denote by $M(B)$ and $M(b)$ the sets of characters of height 0 of $B$ and $b$ respectively. The Alperin-McKay Conjecture then asserts that $|M(B)|=|M(b)|$ (while the McKay Conjecture states that $|M(G)| = |M(N_G(P))|$, where $P \in Syl_p(G)$, and $M(G)$ and $M(N_G(P))$ denote the sets of irreducible characters of $p'$-degree of $G$ and $N_G(P)$ respectively).

In \cite{Isaacs-Navarro}, Isaacs and Navarro predicted that something stronger must happen, namely that this equality can be refined when considering the $p'$-parts of the character degrees. For any $n \in \N$, we write $n=n_pn_{p'}$, with $n_p=p^{\nu(n)}$. For any $1 \leq k \leq p-1$, we define subsets $M_k(B)$ and $M_k(b)$ of $M(B)$ and $M(b)$ respectively by letting $M_k(B)=\{ \chi \in M(B) \, ; \, \chi(1)_{p'} \equiv \pm k \modp \} $ and $ M_k(b)=\{ \varphi \in M(b) \, ; \, \varphi(1)_{p'} \equiv \pm k \modp \}$. We then have the following
\begin{conjecture}\cite[Conjecture B]{Isaacs-Navarro}\label{INConj}
For $1 \leq k \leq p-1$, we have $|M_{ck}(B)|=|M_k(b)|$, where $c=[G \colon N_G(D)]_{p'}$.
\end{conjecture}
Note that Conjecture \ref{INConj} obviously implies the Alperin-McKay Conjecture (by letting $k$ run through $\{1, \, \ldots , \, p-1\}$), but also implies another refinement of the McKay Conjecture; if we let $M_k(G)=\{ \chi \in \Irr(G) \, ; \, \chi(1) \equiv \pm k \modp \}$, then, by considering all blocks of $G$ with defect group $P \in Syl_p(G)$, we obtain $|M_k(G)|=|M_k(N_G(P))|$, since $[G \colon N_G(P)] \equiv 1 \modp$ (see \cite[Conjecture B]{Isaacs-Navarro}).

Isaacs and Navarro proved Conjecture \ref{INConj} whenever $D$ is cyclic, or $G$ is $p$-solvable or sporadic. P. Fong proved it for symmetric groups $S(n)$ in \cite{Fong}, and R. Nath for alternating groups $A(n)$ in \cite{Nath}. In this paper, we prove that Conjecture \ref{INConj} holds in all the covering groups of the symmetric and alternating groups, provided $p$ is odd (Theorem \ref{INC}).

In Section 2, we present the covering groups $S^+(n)$ and $S^-(n)$ and their irreducible characters, first studied by I. Schur in \cite{Schur}, as well as their $p$-blocks. It turns out that the main work to be done is on so-called spin blocks. We also give various results on the degrees of spin characters, generalizing the methods used by Fong in \cite{Fong}. Section 3 is devoted to proving Theorem \ref{Reduction Theorem} which reduces the problem to proving only that Conjecture \ref{INConj} holds for the principal spin block of $S^+(pw)$. This reduction theorem is a refinement of \cite[Theorem 2.2]{Michler-Olsson} that G. O. Michler and J. B. Olsson proved in order to establish that the Alperin-McKay Conjecture holds for covering groups. Finally, the case of the principal spin block of $S^+(pw)$ is treated in Section 4.

\section{Covering groups}\label{part2}
%We keep the notation as in Section~\ref{intro}.

In this section, we introduce the objects and preliminary results we will need about covering groups and their characters. Unless stated otherwise, the following results can be found in \cite{Michler-Olsson}.

\subsection{Covering groups}
For any integer $n \geq 1$, I. Schur has defined (by generators and relations) two central extensions $\hat{S}(n)$ and $\tilde{S}(n)$ of the symmetric group $S(n)$ (see \cite{Schur}, p. 164). We have  $\hat{S}(1) \cong \tilde{S}(1) \cong \Z/2\Z$, and, for $n \geq 2$, there is a nonsplit exact sequence$$1 \longrightarrow \cyc{z} \longrightarrow \hat{S}(n)  \mathop{\longrightarrow}^{\pi} S(n) \longrightarrow 1,$$where $\cyc{z}=Z(\hat{S}(n)) \cong \Z/2\Z$.

Whenever $n \geq 2$, these two extensions are non-isomorphic, except when $n=6$. However, their character tables are the same. Hence, for our purpose, it is sufficient to study one of them. Throughout this paper, we will write $S^+(n)$ for $\hat{S}(n)$.

If $H$ is a subgroup of $S(n)$, we let $H^+= \pi^{-1}(H)$ and $H^-=\pi^{-1}(H \cap A(n))$. In particular, $H^-$ has index 1 or 2 in $H^+$, and $H^+=H^-$ if and only if $H \subset A(n)$. We define $S^-(n)=A(n)^-=A(n)^+$. Hence $S^-(n)$ is a central extension of $A(n)$ of degree 2.

The groups $A(6)$ and $A(7)$ also have one 6-fold cover each, which, together with the above groups, give all the covering groups of $S(n)$ and $A(n)$.

\subsection{Characters, blocks and twisted central product}
From now on, we fix an odd prime $p$. For any $H \leq S(n)$, the irreducible complex characters of $H^{\e}$ fall into two categories; those that have $z$ in their kernel, and which can be identified with those of $H$ (if $\e=1$) or those of $H\cap A(n)$ (if $\e=-1$), and those that don't have $z$ in their kernel. These (faithful) characters are called {\emph{spin characters}}. We denote by $SI(H^{\e})$ the set of spin characters of $H^{\e}$, and we let $SI_0(H^{\e})=SI(H^{\e}) \cap M(H^{\e})$ (with the notation of Section \ref{intro}).

If $B$ is a $p$-block of $H^{\e}$, then, because $p$ is odd, it is known that either $B \cap SI(H^{\e})=\emptyset$ or $B \subset SI(H^{\e})$, in which case we say that $B$ is a {\emph{spin block}} of $H^{\e}$.

\medskip
Any two $\chi, \, \psi \in \Irr(H^{\e})$ are called {\emph{associate}} if $\chi \uparrow^{H^+}= \psi \uparrow^{H^+}$ (if $\e=-1$) or if $\chi \downarrow_{H^-}=\psi \downarrow_{H^-}$ (if $\e=1$). Then each irreducible character of $H^{\e}$ has exactly 1 or 2 associate characters. If $\chi$ is itself its only associate, we say that $\chi$ is {\emph{self-associate}} (written s.a.), we put $\chi^a=\chi$, and let $\sa(\chi)=1$. Otherwise, $\chi$ has a unique associate $\psi \neq \chi$; we say that $\chi$ is {\emph{non-self-associate}} (written n.s.a.), we put $\chi^a= \psi$, and we let $\sa(\chi)=-1$.

If $H^+ \neq H^-$, then $\chi \in \Irr(H^+)$ and $\varphi \in \Irr(H^-)$ are said to {\emph{correspond}} if $\cyc{\chi, \, \varphi \uparrow^{H^+} }_{H^+} \neq 0$. In this case, Clifford's theory implies that $\sa(\chi)=-\sa(\varphi)$.

\medskip

If $H_1, \, H_2, \, \ldots, \, H_k \leq S(n)$ act (non-trivially) on disjoint subsets of $\{1, \, \ldots , \, n\}$, then one can define the {\emph{twisted central product}} $H^+=H_1^+ \hat{\times} \cdots \hat{\times} H_k^+ \leq S^+(n)$ (see \cite{Schur} or \cite{Humphreys}). Then $|H^+| = \frac{1}{2^{k-1}} |H_1^+| . |H_2^+| \ldots |H_k^+| = 2 |H_1| . |H_2| \ldots |H_k|$. Also, one obtains $SI(H^+)$ from the $SI(H_i^+)$'s as follows:

\begin{proposition}\cite[\S28]{Schur}
There is a surjective map
$$\hat{\otimes} \colon \left\{ \begin{array}{ccc} SI(H_1^+) \times \cdots \times SI(H_k^+) & \longrightarrow & SI(H^+) \\ (\chi_1, \, \ldots , \, \chi_k) & \longmapsto & \chi_1 \hat{\otimes} \cdots  \hat{\otimes} \chi_k \end{array} \right.$$which satisfies the following properties.
Suppose $\chi_i, \psi_i \in SI(H_i^+)$ for $1 \leq i \leq k$. Then

\noindent
(i) $\sa(\chi_1 \hat{\otimes} \cdots  \hat{\otimes} \chi_k )=\sa(\chi_1)  \ldots  \sa( \chi_k)$, and $(\chi_1 \hat{\otimes} \cdots  \hat{\otimes} \chi_k  )(1)=2^{\floor{s/2}} \chi_1(1) \ldots  \chi_k(1)$, where $s$ is the number of n.s.a. characters in $\{\chi_1 , \,  \ldots  , \, \chi_k \}$ and $\floor{ \;  }$ denotes integral part.

\noindent
(ii) $\chi_1 \hat{\otimes} \cdots  \hat{\otimes} \chi_k$ and $\psi_1 \hat{\otimes} \cdots  \hat{\otimes} \psi_k $ are associate if and only if $\chi_i$ and $\psi_i$ are associate for all $i$.

\noindent
(iii) $\chi_1 \hat{\otimes} \cdots  \hat{\otimes} \chi_k = \psi_1 \hat{\otimes} \cdots  \hat{\otimes} \psi_k $ if and only if $\chi_i$ and $\psi_i$ are associate for all $i$ and [$\sa(\chi_1)  \ldots  \sa( \chi_k)=1$] or [$\sa(\chi_1)  \ldots  \sa( \chi_k)=-1$ and $| \{ i \, | \, \chi_i \neq \psi_i \} |$ is even].

\end{proposition}

\subsection{Partitions and bar-partitions}

Just as the irreducible characters of $S(n)$ are parametrized by the partitions of $n$, the spin characters of $S^+(n)$ have a combinatorial description. We let $P(n)$ be the set of all partitions of $n$, and $P_0(n)$ be the subset of all partitions in distinct parts, also called {\emph{bar-partitions}}. We write $\la \vdash n $ for $\la \in P(n)$, and $\la \succ n$ for $\la \in P_0(n)$. We also write, in both cases, $|\la|=n$.

\medskip

It is well known that $\Irr(S(n))= \{ \chi_{\la}, \, \la \vdash n \}$. For any $\la \vdash n$, we write $h(\la)$ for the product of all hook-lengths in $\la$. We then have $h(\la)=h_{\la,p}h_{\la,p'}$, where $h_{\la,p}$ (respectively $h_{\la,p'}$) is the product of all hook-lengths divisible by $p$ (prime to $p$ respectively) in $\la$. The Hook-Length Formula then gives $\chi_{\la}(1)=\frac{n!}{h(\la)}$.

If we remove all the hooks of length divisible by $p$ in $\la$, we obtain its $p${\emph{-core}} $\la_{(p)}$. The information on $p$-hooks is stored in the $p${\emph{-quotient}} $\la^{(p)}$ of $\la$. If $n=pw+r$, with $\la_{(p)} \vdash r$, then $\la^{(p)}$ is a $p$-tuple of partitions of $w$, i.e. $\la^{(p)}=(\la^{(0)}, \, \ldots , \, \la^{(p-1)})$ and $|\la^{(0)}| +  \cdots + | \la^{(p-1)}| = w$. The partition $\la$ is uniquely determined by its $p$-core and $p$-quotient. Also, for any integer $k$, there exists a (canonical) bijection between the $kp$-hooks in $\la$ and the $k$-hooks in $\la^{(p)}$ (i.e. in the $\la^{(i)}$'s). 

Finally, the Nakayama Conjecture states that $\chi_{\la}, \, \chi_{\mu} \in \Irr(S(n))$ belong to the same $p$-block if and only if $\la$ and $\mu$ have the same $p$-core.

\medskip
We now present the analogue properties for bar-partitions and spin characters. For any bar-partition $\la=(a_1, \, \ldots , \, a_m)$ of $n$, with $a_1> \cdots  > a_m>0$, we let $m(\la)=m$, and define the {\emph{sign}} of $\la$ by $\sa(\la)=(-1)^{n-m(\la)}$. We then have

\begin{theorem}\cite[\S 41]{Schur}
For each sign $\e \in \{1, \, -1\}$, there is a (canonical) surjective map $f^{\e} \colon SI(S^{\e}(n)) \longrightarrow P_0(n)$ such that:

(i) $\sa(\chi) = \e \sa(f^{\e}(\chi))$ for all $\chi \in SI(S^{\e}(n))$.

(ii) For any $\chi, \, \psi \in SI(S^{\e}(n))$, we have $f^{\e}(\chi)=f^{\e}(\psi)$ if and only if $\chi$ and $\psi$ are associate.

(iii) If $\chi \in SI(S^+(n))$ and $\varphi \in SI(S^-(n))$, then $f^+(\chi)=f^-(\varphi)$ if and only if $\chi$ and $\varphi$ correspond.

\end{theorem}
In particular, each $\la \succ n$ labels one s.a. character $\chi$ or two associate characters $\chi$ and $\chi^a$. Throughout this paper, we will denote by $\cla$ the set of spin characters labeled by $\la$, and write (abusively) $\cla \in SI(S^{\e}(n))$, and $\cla (1)$ for the (common) degree of any spin character in $\cla$. We will also sometimes write $\cla_+$ to emphasize that $\cla \in SI(S^+(n))$ (and $\cla_-$ if $\cla \in SI(S^-(n))$).

\smallskip
For the following results on bars, cores and quotients, we refer to \cite{Olsson-FrobeniusSymbols}. For any odd integer $q$, let $e=(q-1)/2$. We define a $\bar{q}${\emph{-quotient of weight}} $w$ to be any tuple of partitions $(\la^{(0)}, \, \la^{(1)}, \, \ldots , \, \la^{(e)})$ such that $\la^{(0)} \in P_0(w_0)$, $\la^{(i)} \in P(w_i)$ for $1\leq i \leq e$, and $w_0 + w_1 + \cdots + w_e=w$. We define its {\emph{sign}} by $\sa((\la^{(0)}, \, \la^{(1)}, \, \ldots , \, \la^{(e)}))=(-1)^{w-w_0} \sa(\la^{(0)})$.

Now take any bar-partition $\la=(a_1, \, \ldots , \, a_m)$ of $n$ as above. The {\emph{bars}} in $\la$ can be read in the {\emph{shifted Young diagram}} $S(\la)$ of $\la$. This is obtained from the usual Young diagram of $\la$ by shifting the $i$-th row $i-1$ positions to the right. The $j$-th node in the $i$-th row is called the $(i,j)$-node, and correspond to the {\emph{bar}} $B_{ij}$. The {\emph{bar-lengths}} in the $i$-th row are obtained by writing (from left to right in $S(\la)$) the elements of the following set in decreasing order: $\{1, \, 2, \, \ldots , \,  a_i \} \cup \{a_i + a_j \, | \, j > i \} \setminus \{ a_i - a_j \, | \, j>i \}$. The bars are of three types:
\begin{itemize}
\item{}
Type 1. These are bars $B_{ij}$ with $i+j \geq m+2$ (i.e. in the right part of $S(\la)$). They are ordinary hooks in $S(\la)$, and their lengths are the elements of $\{1, \, 2, \, \ldots , \,  a_i -1\} \setminus \{ a_i - a_j \, | \, j>i \}$.
\item{}
Type 2. These are bars $B_{ij}$ with $i+j = m+1$ (in particular, the corresponding nodes all belong to the same column of $S(\la)$). Their length is precisely $a_i$, and the bar is all of the $i$-th row of $S(\la)$.
\item
Type 3. The lengths $\{a_i + a_j \, | \, j > i \} $ correspond to bars $B_{ij}$ with $i+j \leq m$. The bar consists of the $i$-th row together with the $j$-th row of $S(\la)$.
\end{itemize}
Bars of type 1 and 2 are called {\emph{unmixed}}, while those of type 3 are called {\emph{mixed}}. The unmixed bars in $\la$ correspond exactly to the hooks in the partition $\la^*$, which admits as a $\beta$-set the set of parts of $\la$.

\smallskip
For any $\la \succ n$, we write $\h(\la)$ for the product of all bar-lengths in $\la$. We then have $\h(\la)=\h_{\la,p}\h_{\la,p'}$, where $\h_{\la,p}$ (respectively $\h_{\la,p'}$) is the product of all bar-lengths divisible by $p$ (prime to $p$ respectively) in $\la$. We then have the following analogue of the Hook-Length Formula (proved by A. O. Morris, \cite[Theorem 1]{Morris})
$$\cla(1)= 2^{\floor{(n-m(\la))/2}}\frac{n!}{\h(\la)}.$$

\smallskip
If we remove all the bars of length divisible by $p$ in $\la$, we obtain its $\p${\emph{-core}} $\la_{(\p)}$ (which is still a bar-partition), and its $\p${\emph{-quotient}} $\la^{(\p)}$. If $n=pw+r$, with $\la_{(\p)} \succ r$, then $\la^{(\p)}$ is a $\p$-quotient of weight $w$ in the sense defined above. The bar-partition $\la$ is uniquely determined by its $\p$-core and $\p$-quotient. Also, for any integer $k$, there exists a canonical bijection between the set of $kp$-bars in $\la$ and the set of $k$-bars in $\la^{(\p)}$ (where a $k$-bar in $\la^{(\p)}= (\la^{(0)}, \, \la^{(1)}, \, \ldots , \, \la^{((p-1)/2)})$ is a $k$-bar in $\la^{(0)}$ or a $k$-hook in one of $\la^{(1)}, \, \ldots , \, \la^{((p-1)/2)}$). 

\smallskip
The distribution of the spin characters of $S^+(n)$ into spin blocks was first conjectured for $p$ odd by Morris. It was first proved by J. F. Humphreys in \cite{Humphreys-Blocks}, then differently by M. Cabanes, who also determined the structure of the defect groups of spin blocks (see \cite{Cabanes}).

\begin{proposition}
Let $\chi, \, \psi \in SI(S^{\e}(n))$ and $p$ be an odd prime. Then $\chi$ is of $p$-defect 0 if and only if $f^{\e}(\chi)$ is a $\p$-core. If $f^{\e}(\chi)$ is not a $\p$-core, then $\chi$ and $\psi$ belong to the same $p$-block if and only if $f^{\e}(\chi)_{(\p)}=f^{\e}(\psi)_{(\p)}$.
\end{proposition}
One can therefore define the $\p$-core of a spin block $B$ and its {\emph{weight}} $w(B)$, as well as its {\emph{sign}} $\da(B)=\sa( f^{\e}(\chi)_{(\p)})$ (for any $\chi \in B$). We then have
\begin{proposition}\cite{Cabanes}
If $B$ is a spin block of $S^{\e}(n)$ of weight $w$, then a defect group $X$ of $B$ is a Sylow $p$-subgroup of $S^{\e}(pw)$.
\end{proposition}

\subsection{Removal of $p$-bars}

The following result is the bar-analogue of \cite[Lemma 3.2]{Fong}; it describes how the removal of $p$-bars affects the product of $p'$-bar-lengths.
\begin{proposition}\label{Fong's Lemma}
Suppose $\lambda \succ n$ has $\bar{p}$-core $\lambda_{(\bar{p})}$. Then
$$\bar{h}_{\lambda, p'} \equiv \pm 2^{-a(\lambda)}  \bar{h}_{\lambda_{(\bar{p})}, p'} =  \pm 2^{-a(\lambda)}  \bar{h}(\lambda_{(\bar{p})}) \; \; (\mbox{mod} \; p),$$
where $a(\lambda)$ is the number of $p$-bars of type 3 to remove from $\lambda$ to get $\lambda_{(\bar{p})}$.

\end{proposition}
\begin{proof}
Let $B_{ij}$ be a $p$-bar in $\lambda$ and $\la - B_{ij}$ be the bar-partition obtained from $\la$ by removing $B_{ij}$. We distinguish two cases, depending on whether $B_{ij}$ is unmixed or mixed.

\smallskip
First suppose that $B_{ij}$ is unmixed (i.e. $i+j > m(\la)$). We start by examining the unmixed $p'$-bars in $\la$ and $\la-B_{ij}$. These correspond, in the notation above, to the $p'$-hooks in $\la^*$ and $(\la-B_{ij})^*$ respectively (considering  $\la$ and $\la-B_{ij}$ as $\beta$-sets). The set of parts of $\la$ is $X= \{ a_1, \, \ldots , \, a_m \}$, and the set of non-zero parts of $\la-B_{ij}$ is $Y=\{a_1, \, \ldots , \, a_{i-1}, \, a_i -p, \, a_{i+1}, \, \ldots , \, a_m\}$ (or $Y=\{a_1, \, \ldots , \, a_{i-1}, \, a_{i+1}, \, \ldots , \, a_m\}$ if $a_i=p$). The $p'$-hooks in $\la^*$ (resp. $(\la-B_{ij})^*$) therefore correspond to pairs $(x, \, y)$, with $0 \leq x <y$, $( y-x, \, p)=1$, and $x\not \in X$, $y \in X$ (resp. $x \not \in Y$, $y \in Y$).

If $B_{ij}$ is of type 1 (i.e. $i + j > m(\la)+1$), then $a_i - p >0$, so that $|Y|=|X|$ and $(\la-B_{ij})^*= \la^* - h$ for some $p$-hook $h$ in $\la$. In this case, we are thus exactly in the same context as \cite[Lemma 3.2]{Fong}, and we get $h_{\la^*,p'} \equiv - h_{\la^*-h, p'}=-h_{(\la-B_{ij})^*,p'} \; (\mbox{mod} \; p)$.

If, on the other hand, $B_{ij}$ is of type 2 (i.e. $i+j = m(\la)+1$), then $a_i -p =0$, and $Y=X \setminus \{p\}$. Note that, in this case, $Y$ is not a $\beta$-set for a partition of $|\la^*|-p$, while $Y \cup \{0\}$ is. The $p'$-hooks in $(\la - B_{ij})^*$ correspond to either pairs $(x, \, y)$ with $y \neq a_i$, which also correspond to $p'$-hooks in $\la^*$, or to pairs $(p, \, y)$, with $y >p$ and $y \in X$. These new hooks have lengths $(a_1-p), \, \ldots, \, (a_{i-1}-p)$. Finally, some hooks have disappeared: those corresponding to pairs $(x, \, p)$ with $x<p$ and $x \not \in X$. These have lengths $(p-x)$, for $0 \leq x <p$ and $x \not \in \{ a_{i+1}, \, \ldots , \, a_m \}$.

We now turn to the mixed $p'$-bars in $\la$ and $\la-B_{ij}$. Suppose first that $B_{ij}$ is of type 1. Then $m(\la)=m(\la-B_{ij})$. Suppose that
$$a_1 > \cdots >a_{i-1} > a_{i+1}>\cdots >a_k > a_i-p>a_{k+1}>\cdots >a_m.$$
To prove the result, we can simply ignore the bar-lengths which are common to $\la$ and $\la-B_{ij}$. The mixed bars which disappear when going from $\la$ to $\la-B_{ij}$ have lengths $$(a_1+a_i), \, (a_2+a_i), \, \ldots, \, (a_{i-1}+a_i)$$
$$\mbox{and} \;  (a_i + a_{i+1}), \, (a_i + a_{i+2}), \, \ldots , \, (a_i+a_m).$$
The mixed bars which appear have lengths$$(a_1+a_i-p), \, (a_2+a_i-p), \, \ldots, \, (a_{i-1}+a_i-p),$$
$$ (a_{i+1}+a_i-p),  \, \ldots , \, (a_k + a_i -p) \; \mbox{and} \;  (a_i -p + a_{k+1}),  \, \ldots , \, (a_i-p+a_m).$$
If we then just consider the lengths not divisible by $p$, it is easy to see that we can pair the bars disappearing with those appearing. The pairs are of the form $(b,b')$, where $b$ is a bar in $\la$ and $b'$ is a bar in $\la-B_{ij}$, and $|b'|=|b|-p$. We thus get, in this case,$$\prod_{b \, \mbox{mixed} \, p'\mbox{-bar in} \, \la} |b| \equiv \prod_{b' \, \mbox{mixed} \, p'\mbox{-bar in} \, \la-B_{ij}} |b'| \; \; (\mbox{mod} \; p).$$ 
Together with the equality obtained above for unmixed $p'$-bars, we obtain that, if $B_{ij}$ is a $p$-bar of type 1 in $\lambda$, then $\bar{h}_{\la, p'} \equiv \pm \bar{h}_{\la-B_{ij}, p'} \; (\mbox{mod} \; p)$.

Now suppose that $B_{ij}$ is of type 2, i.e. $a_i=p$. Then the mixed bars which disappear when going from $\la$ to $\la-B_{ij}$ have lengths $(a_1+p), \, (a_2+p), \, \ldots, \, (a_{i-1}+p)$ (call these $A$) and $(p+ a_{i+1}), \, (p + a_{i+2}), \, \ldots , \, (p+a_m)$ (call these $B$), while no new mixed bar appears.

The bars disappearing in $A$ are compensated for by the hooks appearing in $(\la-B_{ij})^*$ in the study of unmixed bars above (since $a_u+p \equiv a_u -p \, \modp$ for all $1 \leq u \leq i-1$, the $p'$-parts are congruent mod $p$ when these are not divisible by $p$).

On the other hand, since $0<a_m < \cdots < a_{i+1} <a_i =p$, all the bar-lengths in $B$ are coprime to $p$, and their product is
$$(p+ a_{i+1})(p + a_{i+2}) \ldots  (p+a_m) \equiv a_{i+1} a_{i+2} \ldots a_m \; \; (\mbox{mod} \; p).$$
Now the hooks disappearing in the above discussion of unmixed bars all have length prime to $p$ except one (corresponding to $x=0$). The product of the lengths prime to $p$ is thus
$$\prod_{0<x<p, \, x \not \in \{a_i+1, \, \ldots , \, a_m \} } (p-x) \equiv  (-1)^{p-1-m+i} \prod_{0<x<p, \, x \not \in \{a_i+1, \, \ldots , \, a_m \} } x \; \; (\mbox{mod} \; p).$$Hence the product of the $p'$-hook-lengths disappearing and the $p'$-bar-lengths in $B$ is congruent (mod $p$) to$$(-1)^{p-1-m+i} \prod_{0<y<p} y=(-1)^{p-1-m+i} (p-1)! \equiv (-1)^{p-m+i} \modp$$(by Wilson's Theorem). Finally, we obtain that, if $B_{ij}$ is a $p$-bar of type 2 in $\lambda$, then $\bar{h}_{\la, p'} \equiv (-1)^{p-m+i} \bar{h}_{\la-B_{ij}, p'} \; (\mbox{mod} \; p)$.

\smallskip
We now suppose that $B_{ij}$ is a $p$-bar of type 3 in $\lambda$, i.e. $i<j$, $a_i > a_j$ and $a_i + a_j=p$. The set of parts of $\la$ is $X=\{a_1, \, \ldots a_m \}$ and the set of parts of $\la-B_{ij}$ is $Y=\{ a_1 , \, \ldots , \, a_{i-1}, \, a_{i+1}, \, \ldots , \, a_{j-1}, \, a_{j+1}, \, \ldots , \, a_m\}$.
Ignoring as before the bars which are common to $\la$ and $\la - B_{ij}$, we see that the unmixed bars which disappear from $\la$ to $\la - B_{ij}$ have lengths
$$(a_i-x) \;  \;  \; (0 \leq x <a_i, \; x \not \in \{a_m, \, \ldots , \, a_{i+1} \}) $$
$$ \mbox{and} \; \; (a_j-x) \;  \;  \; (0 \leq x <a_j, \; x \not \in \{a_m, \, \ldots , \, a_{j+1} \}),$$
while those appearing have lengths
$$(a_1-a_i), \, \ldots , \, (a_{i-1}-a_i), \; (a_1 - a_j), \, \ldots , \, (a_{i-1} - a_j), \; \mbox{and} \; (a_{i+1}-a_j), \, \ldots , \, (a_{j-1}-a_j).$$
On the other hand, there is no mixed bar appearing, while the mixed bars disappearing have lengths
$$(a_1+a_i), \, \ldots , \, (a_{i-1}+a_i) \; \; \; \; \; \mbox{(rows 1, ... $i-1$, column $i$)},$$
$$(a_i+a_{i+1}), \, \ldots, \, (a_i+a_{j-1}), \, (a_i+a_j), \, \ldots , \, (a_i + a_m) \; \; \; \; \mbox{(row $i$)},$$
$$(a_1+a_j), \, \ldots , \, (a_{i-1}+a_j) \; \; \; \; \; \mbox{(rows 1, ... $i-1$, column $j$)},$$
$$(a_{i+1}+a_j), \, \ldots , \, (a_{j-1}+a_j) \; \; \; \; \; \mbox{(rows $i+1$, ... $j-1$, column $j$)},$$
$$\mbox{and} \; (a_j+a_{j+1}), \, \ldots, \, (a_j + a_m) \; \; \; \; \mbox{(row $j$)}.$$
Now, since $a_i+a_j=p$, we have, for any $1 \leq k \leq m$,
$$a_k-a_i \equiv a_k+a_j \modp \; \; \mbox{and} \; \; a_k+a_i \equiv a_k-a_j \modp.$$
In particular, $a_k-a_i$ (resp. $a_k+a_i$) is coprime to $p$ if and only if $a_k+a_j$ (resp. $a_k-a_j$) is coprime to $p$, and, in that case,
$$(a_k\pm a_i)_{p'} = a_k \pm a_i \equiv a_k \mp a_j \equiv (a_k \mp a_j)_{p'} \modp .$$
We thus have the following compensations between the appearing unmixed bars and the appearing mixed bars:
$$(a_1-a_i), \, \ldots , \, (a_{i-1}-a_i) \; \; \longleftrightarrow \; \; (a_1+a_j), \, \ldots , \, (a_{i-1}+a_j) $$
$$(a_1 - a_j), \, \ldots , \, (a_{i-1} - a_j) \; \; \longleftrightarrow \; \;  (a_1+a_i), \, \ldots , \, (a_{i-1}+a_i)$$
$$ \mbox{and} \; (a_{i+1}-a_j), \, \ldots , \, (a_{j-1}-a_j)  \; \; \longleftrightarrow \; \;  (a_i+a_{i+1}), \, \ldots, \, (a_i+a_{j-1}).$$
This accounts for all the appearing (unmixed) bars, and we're left exactly with the following disappearing bar-lengths:
$$(a_i - x ) \; \; (0 \leq x <a_i, \; x \not \in \{a_m, \, \ldots , \, a_{i+1} \}) \; \; \mbox{unmixed of type 1},$$
$$(a_j-x) \;  \;  (0 \leq x <a_j, \; x \not \in \{a_m, \, \ldots , \, a_{j+1} \})  \; \; \mbox{unmixed of type 2},$$
$$(a_{i+1}+a_j), \, \ldots , \, (a_{j-1}+a_j) , \,  (a_j+a_{j+1}), \, \ldots, \, (a_j + a_m) \; \; \mbox{mixed of type 1},$$
$$(a_i+a_{j+1}), \, \ldots , \, (a_i + a_m)\; \; \mbox{mixed of type 2},$$
and $(a_i + a_j)=p$ which can thus be ignored.

Now, for any $i+1 \leq k \leq m$, $a_j+a_k = p-a_i+a_k \equiv -(a_i-a_k) \modp$, and, for $j+1 \leq k \leq m$, $a_i+a_k \equiv -(a_j-a_k) \modp$. Hence, taking the product, we obtain (modulo $p$):
$$\prod_{0 \leq x < a_i, \, x \neq a_j} (a_i-x)=\frac{a_i!}{a_i-a_j} \; \; (\mbox{type 1}) \; \; \mbox{and} \; \;     \prod_{0 \leq x < a_j} (a_j-x)=a_j! \; \; (\mbox{type 2}).$$
Now $a_j!=1.2\ldots a_j=(-1)^{a_j}(-1)\ldots(-a_j) \equiv (-1)^{a_j}(p-1)\ldots(p-a_j) \modp$, so that $a_j!\equiv (-1)^{a_j}(p-1)\ldots(a_i+1)a_i \modp$.  We thus have, disappearing,
$$\displaystyle \pm \frac{a_i a_i! (a_i+1) \ldots (p-1)}{a_i-a_j} \equiv \pm \frac{a_i}{a_i-a_j} (p-1)! \equiv \mp  \frac{a_i}{a_i-a_j} \modp$$
(this last equality being true by Wilson's Theorem).

Finally, $a_i-a_j=a_i-(p-a_i) \equiv -2 a_i \modp$, yielding a total of $\pm 2^{-1} \modp$ disappearing (since, $p$ being odd, 2 is invertible$\modp$, and $a_i <p$ so that we can simplify by $a_i$). We thus get that, if $B_{ij}$ is a $p$-bar of type 3 in $\lambda$, then $\bar{h}_{\la, p'} \equiv \pm \frac{1}{2} \bar{h}_{\la-B_{ij}, p'} \; (\mbox{mod} \; p)$.

\smallskip
Iterating the above results on all the $p$-bars to remove from $\la$ to get to its $\bar{p}$-core $\la_{(\bar{p})}$, we finally obtain the desired equality, writing $a(\lambda)$ for the number of $p$-bars of type 3 to remove:
$$\bar{h}_{\lambda, p'} \equiv \pm 2^{-a(\lambda)}  \bar{h}_{\lambda_{(\bar{p})}, p'} =  \pm 2^{-a(\lambda)}  \bar{h}(\lambda_{(\bar{p})}) \; \; (\mbox{mod} \; p)$$
(since all the bars in $\la_{(\bar{p})}$ have length coprime to $p$).

\end{proof}

\subsection{$\p$-core tower, $\p$-quotient tower and characters of $p'$-degree}

In this section, we want to obtain an expression for the (value modulo $p$ of the) $p'$-part of the degree of a spin character. We start by describing the $\p$-core tower of a bar-partition, introduced by Olsson in \cite{Olsson-FrobeniusSymbols}.

Take any $\la \succ n$. the $\p$-core tower of $\la$ has rows $R_0^{\la}, \, R_1^{\la}, \, R_2^{\la}, \, \ldots$, where the $i$-th row $R_i^{\la}$ contains one $\p$-core and $(p^i-1)/2$ $p$-cores (in particular, one can consider $R_i^{\la}$ as a $\p^i$-quotient). We have $R_0^{\la} = \{ \la_{(\p)} \}$ (the $\p$-core of $\la$). If the $\p$-quotient of $\la$ is $\la^{(\p)}=(\la^{(0)}, \, \la^{(1)}, \, \ldots , \, \la^{(e)}) $ (where $e=(p-1)/2$), then $R_1^{\la}= \{ \la_{ (\p)}^{(0)}, \, \la_{ (p)}^{(1)}, \, \ldots , \, \la_{ (p)}^{(e)} \}$. Writing $\la^{(0)(\p)}=(\la^{(0,0)}, \, \la^{(0,1)}, \, \ldots , \, \la^{(0,e)})$ the $\p$-quotient of $\la^{(0)}$ and $\la^{(i)(p)}= (\la^{(i,1)}, \, \la^{(i,2)}, \, \ldots , \, \la^{(i,p)})$ the $p$-quotient of $\la^{(i)}$ ($1 \leq i \leq e$), and taking cores, we let $$R_2^{\la}= \{ \la^{(0,0)}_{ (\p)}, \, \la^{(0,1)}_{ (p)}, \, \ldots , \, \la^{(0,e)}_{ (p)} , \, \la^{(1,1)}_{ (p)}, \,  \ldots , \, \la^{(1,p)}_{ (p)}, \, \la^{(2,1)}_{ (p)} , \, \ldots , \, \la^{(e,p)}_{ (p)} \}.$$Continuing in this way, we obtain the $\p$-core tower of $\la$. We define the $\p$-quotient tower of $\la$ in a similar fashion: it has rows $Q_0^{\la}, \, Q_1^{\la}, \, Q_2^{\la}, \, \ldots$, where the $i$-th row $Q_i^{\la}$ contains one $\p$-quotient and $(p^i-1)/2$ $p$-quotients (in particular, $Q_i^{\la}$ can be seen as a $\p^{i+1}$-quotient). With the above notation, we have $Q_0^{\la}=\{ \la^{(\p)} \}$, $Q_1^{\la}= \{ \la^{(0) (\p)}, \, \la^{(1) (p)}, \, \ldots , \, \la^{(e) (p)} \}$ and $$ Q_2^{\la}= \{ \la^{(0,0) (\p)}, \, \la^{(0,1) (p)}, \, \ldots , \, \la^{(0,e )(p)} , \, \la^{(1,1 )(p)}, \,  \ldots , \, \la^{(1,p) (p)}, \, \la^{(2,1 )(p)} , \, \ldots , \, \la^{(e,p) (p)} \}.$$

\smallskip
The following result will be useful later. 
\begin{lemma}\label{tower-sign}
If $\la \succ n$ has $\p$-core tower $(R_0^{\la}, \, R_1^{\la}, \, \ldots , \, R_m^{\la})$, then $\sa(\la)= \prod_{i=0}^{m}\sa(R_i^{\la})$.
\end{lemma}
\begin{proof}
We have $\sa(\la)=\sa(\la_{(\p)}) \sa(\la^{(\p)})$, and $\sa(\la_{(\p)})=\sa(R_0^{\la})$.

Also, $\sa(\la^{(\p)})=\sa(\la^{(0)}) (-1)^{\sum_{i\geq1} | \la^{(i)} |}$, and$$\sa(\la^{(0)})=\sa(\la^{(0)}_{(\p)})\sa(\la^{(0)(\p)})=\sa(Q_0^{\la})=\sa(\la^{(0)}_{(\p)})\sa(\la^{(0,0)}) (-1)^{\sum_{j \geq 1} | \la^{(0,j)} |}.$$
Now $\sa(R_1^{\la})=\sa(\la^{(0)(\p)}) (-1)^{\sum_{i\geq1} | \la^{(i)}_{(\p)} |}$ and $\sa(Q_1^{\la})=\sa(\la^{(0,0)}) (-1)^{\sum_{i\geq0, j \geq 1 } | \la^{(i,j)} |}$, so that$$\begin{array}{rl}\sa(R_1^{\la})\sa(Q_1^{\la}) & =\sa(\la^{(0)(\p)}) \sa(\la^{(0,0)}) (-1)^{\sum_{i\geq1} | \la^{(i)}_{(\p)} |+ \sum_{i\geq0, j \geq 1 } | \la^{(i,j)} |} \\ & = \sa(\la^{(0)})  (-1)^{\sum_{i\geq1} | \la^{(i)}_{(\p)} |+ \sum_{i, j \geq1 } | \la^{(i,j)} |}. \end{array}$$
However, for each $i \geq 1$, we have $|\la^{(i)}_{(\p)} |+ \sum_{ j \geq1 } | \la^{(i,j)} | \equiv |\la^{(i)}_{(\p)} |+ p \sum_{ j \geq1 } | \la^{(i,j)} | \modtwo$ (since $p$ is odd), and $| \la^{(i)}_{(\p)} |+ p \sum_{ j \geq1 } | \la^{(i,j)} | =| \la^{(i)} |$. We therefore get
$$\sa(R_1^{\la})\sa(Q_1^{\la})  \sa(\la^{(0)})  (-1)^{\sum_{i\geq1} | \la^{(i)} |}=\sa(\la^{(\p)}) =\sa(Q_0^{\la}).$$
Finally, we have $\sa(\la)=\sa(R_0^{\la})\sa(Q_0^{\la})$, and $\sa(Q_0^{\la})=\sa(R_1^{\la})\sa(Q_1^{\la}) $, whence $\sa(\la)=\sa(R_0^{\la})\sa(R_1^{\la})\sa(Q_1^{\la}) $. Iterating this process, we deduce the result.

\end{proof}

Now, writing $\beta_i(\la)$ for the sum of the cardinalities of the partitions in $R_i^{\la}$, one shows easily that $| \la |= \sum_{i \geq 0} \beta_i(\la) p^i$ (see \cite{Olsson-FrobeniusSymbols}). Also, one gets the following bar-analogue of \cite[Proposition 1.1]{Fong}:

\begin{proposition}\cite[Proposition (3.1)]{Olsson-FrobeniusSymbols}\label{height}
In the above notation, 
$$\nu_p( \bar{h}(\la) )= \displaystyle \frac{n-\sum_{i \geq 0} \beta_i(\la)}{p-1}.$$
In particular, $\cla$ has $p'$-degree if and only if $\sum_{i \geq 0} \beta_i(\la) p^i$ is the $p$-adic decomposition of $n$.
\end{proposition}

Let $n= \sum_{i=0}^k t_i p^i$ be the $p$-adic decomposition of $n$. For each $0 \leq i \leq k$, let $e_i =(p^i-1)/2$, and write $R_i^{\la}=\{ \mu_i^{(0)}, \, \mu_i^{(1)}, \, \ldots , \, \mu_{i}^{(e_i)} \}$ and $Q_i^{\la}= \{ \la_i^{(0)}, \, \la_i^{(1)}, \, \ldots , \, \la_i^{(e_{i+1})} \}$. Note that $Q_k^{\la}= \{ \emptyset, \, \ldots , \, \emptyset \}$.

\noindent
We let $\h(R_i^{\la})=  \h(\mu_i^{(0)}) \prod_{j=1}^{e_i} h( \mu_i^{(j)})$, and $\h(Q_i^{\la})=  \h(\la_i^{(0)}) \prod_{j=1}^{e_{i+1}} h( \la_i^{(j)})$, and we let $m_i^{(0)}=m(\mu_i^{(0)})$ and $\beta_i=\beta_i(\la)$.

\begin{proposition}\label{degspin}
With the above notation, we have, for any $\la \succ n$,
$$\frac{|S^+(n)|_{p'}}{\cla (1)_{p'}} \equiv \pm \frac{2}{2^{\floor{\frac{S}{2}}} }\prod_{i=0}^k \frac{1}{2^{\floor{ ( \beta_i-m_i^{(0)})/2} }} \h (R_i^{\la}) \modp ,$$
where $S = | \{ 0 \leq i \leq k \, ; \;  \beta_i - m_i^{(0)} \; \mbox{odd} \} |$.
\end{proposition}

\begin{proof}
We have
$$\frac{|S^+(n)|_{p'}}{\cla (1)_{p'}} =  \frac{2}{2^{\lfloor(n-m(\la))/2\rfloor}} \h (\la)_{p'}.$$
Now $\h (\la)_{p'}= \h_{\la,p'} (\h_{\la,p})_{p'}=\h_{\la,p'} (\h (Q_0^{\la})_{p'})$ (since there is a bijection between the set of bars divisible by $p$ in $\la$ and the set of bars in the quotient $Q_0^{\la}$).

\noindent
By Proposition \ref{Fong's Lemma}, we have $\bar{h}_{\lambda, p'} \equiv \pm 2^{-a(\lambda)}  \bar{h}(\lambda_{(\bar{p})}) \equiv \pm 2^{-a(\lambda)} \h(R_0^{\la}) \modp$. Also,
$$\begin{array}{rcl} \h (Q_0^{\la})_{p'} & = & \h_{Q_0^{\la},p'}  ( \h_{Q_0^{\la},p})_{p'} \\ & = & \h_{\la_0^{(0)},p'} h_{\la_0^{(1)},p'} \ldots h_{\la_0^{(e_1)},p'}  \h (Q_1^{\la})_{p'} \\  & \equiv & \pm 2^{-a(\la_0^{(0)})}  \h_{\mu_1^{(0)}} h_{\mu_1^{(1)}} \ldots h_{\mu_1^{(e_1)}}  \h (Q_1^{\la})_{p'} \modp , \end{array}$$
this last equality holding by \ref{Fong's Lemma} (applied to $\la_0^{(0)}$) and by \cite[Lemma 3.2]{Fong} (applied to $\la_0^{(1)}, \, \ldots , \, \la_0^{(e_1)}$). We thus get $\h (Q_0^{\la})_{p'} \equiv  \pm 2^{-a(\la_0^{(0)})}   \h(R_1^{\la})  \h (Q_1^{\la})_{p'}  \modp$, and
$$\h (\la)_{p'} \equiv \pm 2^{-a(\lambda)-a(\la_0^{(0)})}   \h(R_0^{\la})  \h(R_1^{\la})  \h (Q_1^{\la})_{p'}  \modp.$$
Iterating this, until we get to $Q_k^{\la}= \{ \emptyset, \, \ldots , \, \emptyset \}$, we obtain
$$\h (\la)_{p'} \equiv \pm 2^{-\lfloor a(\la)+a(\la_0^{(0)})+a(\la_1^{(0)})+  \cdots + a(\la_{k-1}^{(0)})\rfloor} \prod_{i=0}^k   \h (R_i^{\la}) \modp.$$
On the other hand, repeated use of \cite[Corollary 2.6]{Olsson-FrobeniusSymbols} yields
$$\begin{array}{rcl} m(\la) & = & m_0^{(0)} + m(\la_0^{(0)}) + 2 a(\la) \\   & = & m_0^{(0)} +  m_1^{(0)} + m(\la_1^{(0)}) + 2 a(\la) + 2 a(\la_0^{(0)})  \\ & = & (\cdots) \\ & = &  m_0^{(0)} +  m_1^{(0)} + \cdots + m_k^{(0)} + 2 (a(\la) + a(\la_0^{(0)})  + \cdots + a(\la_{k-1}^{(0)})), \end{array}$$
so that$$\begin{array}{rcl} \lfloor\frac{(n-m(\la))}{2}\rfloor & = & \lfloor \frac{(n- m_1^{(0)} - \cdots - m_k^{(0)})}{2} - (a(\la) + a(\la_0^{(0)})  + \cdots + a(\la_{k-1}^{(0)})) \rfloor \\ & = & \lfloor \frac{(n- (m_1^{(0)} + \cdots + m_k^{(0)})}{2}\rfloor  - (a(\la) + a(\la_0^{(0)})  + \cdots + a(\la_{k-1}^{(0)})) . \end{array}$$
Together with the expression we obtained for $\h (\la)_{p'} $, this gives
$$\frac{|S^+(n)|_{p'}}{\cla (1)_{p'}} \equiv \pm \frac{2}{2^{\lfloor(n-m_0^{(0)}- \cdots - m_k^{(0)})/2\rfloor} }\prod_{i=0}^k \h (R_i^{\la}) \modp.$$
Now recall that $n=\sum_{i=0}^k \beta_i p^i$. Also, for any $1 \leq i \leq k$, we have
$$\begin{array}{rcl} \floor{ \frac{\beta_i p^i - m_i^{(0)}}{2} }  & = & \floor{ \frac{\beta_i (p^i-1) }{2} +  \frac{\beta_i  - m_i^{(0)}}{2}   } \\   & = & \floor{ \frac{(p-1) \beta_i (1 + p + \cdots + p^{i-1}) }{2} +  \frac{\beta_i  - m_i^{(0)}}{2}   }  \\ & = &  \frac{(p-1)}{2}  \beta_i (1 + p + \cdots + p^{i-1})  +  \floor{ \frac{\beta_i  - m_i^{(0)}}{2}   } ,
\end{array}$$
and
$$2^{\frac{(p-1)}{2}  \beta_i (1 + p + \cdots + p^{i-1}) } = (2^{\frac{(p-1)}{2}})^{ \beta_i (1 + \cdots + p^{i-1}) } \equiv (-1)^{ \beta_i (1  + \cdots + p^{i-1}) } \equiv \pm 1 \modp.$$
Hence
$$\begin{array}{rcl} 2^{ \floor{ \frac{n-(m_0^{(0)}+ \cdots + m_k^{(0)})}{2} }} & = & 2^{ \floor{ \frac{ \sum_{i=0}^k (\beta_i p^i - m_i^{(0)}) }{2} } } \\  & = & 2^{ \floor{ \sum_{i=1}^k \frac{p-1}{2} \beta_i (1 + \cdots + p^{i-1}) + \sum_{i=0}^k \frac{\beta_i - m_i^{(0)}}{2} } } \\    & = & 2^{  \sum_{i=1}^k \frac{p-1}{2} \beta_i (1 + \cdots + p^{i-1}) + \floor{ \sum_{i=0}^k \frac{\beta_i - m_i^{(0)}}{2} } }  \\ & \equiv & \pm 2^{ \floor{ \sum_{i=0}^k \frac{\beta_i - m_i^{(0)}}{2} } } \modp.
\end{array}$$
Now
$$\begin{array}{rcl}    \floor{ \sum_{i=0}^k \frac{\beta_i - m_i^{(0)}}{2} } & = &  \floor{ \sum_{i=0 , \; \beta_i - m_i^{(0)} \; \mbox{\tiny{even}}}^k \frac{\beta_i - m_i^{(0)}}{2}  +  \sum_{i=0 , \;  \beta_i - m_i^{(0)} \; \mbox{\tiny{odd }}}^k \frac{\beta_i - m_i^{(0)}}{2} } \\
& = &   \sum_{i=0 , \;  \beta_i - m_i^{(0)} \; \mbox{\tiny{even}}}^k  \floor{ \frac{\beta_i - m_i^{(0)}}{2} }  +  \floor{  \sum_{i=0 , \;  \beta_i - m_i^{(0)} \; \mbox{\tiny{odd }}}^k \frac{\beta_i - m_i^{(0)}}{2} }

\end{array}$$
and we have $\floor{  \sum_{i=0 , \;  \beta_i - m_i^{(0)} \; \mbox{\tiny{odd }}}^k \frac{\beta_i - m_i^{(0)}}{2} } =  \floor{ \frac{S}{2} } +  \sum_{i=0 , \;  \beta_i - m_i^{(0)} \; \mbox{\tiny{odd }}}^k  \floor{  \frac{\beta_i - m_i^{(0)}}{2} } $, where $S = | \{ 0 \leq i \leq k \, ; \;  \beta_i - m_i^{(0)} \; \mbox{odd} \} |$. We finally obtain
$$\frac{|S^+(n)|_{p'}}{\cla (1)_{p'}} \equiv \pm \frac{2}{2^{\floor{\frac{S}{2}}} }\prod_{i=0}^k \frac{1}{2^{\floor{ ( \beta_i-m_i^{(0)})/2} }} \h (R_i^{\la}) \modp .$$

\end{proof}

\section{Reduction Theorem}\label{part3}

In this section, we show that, in order to prove Conjecture \ref{INConj} for any spin block $B$ of $S^{\e}(n)$ of positive weight $w$, it is enough to prove it for the principal spin block of $S^+(pw)$ (i.e. that with empty $\p$-core). Our main tool to navigate between $S^+(n)$ and $S^-(n)$ is the strong duality that exists between their spin blocks.

\smallskip
Let $H \leq S(n)$. A block of $H^{\e}$ is called {\emph{proper}} if it contains both a s.a. character and a n.s.a. character. By \cite[2.1]{olsson-blocks}, any spin block of $S^{\e}(n)$ of positive weight is proper. Now, if $B$ is a proper block of $H^{\e}$, then $H^{\e} \neq H^{-\e}$, and there exists a unique block $B^*$ of $H^{-\e}$ covering $B$ (if $\e=-1$) or covered by $B$ (if $\e=1$), and $B^*$ is also proper. We say that $B$ and $B^*$ are {\emph{(dual) corresponding blocks}}. Finally, if $B$ is proper, then it follows that $B$ consists of s.a. characters and pairs of n.s.a. characters. In particular, we can still write (abusively) $\cla \in B$ or $\cla \in M(B)$. Also, for any sign $\varepsilon$, if $\cla_{\varepsilon} \in SI(S^{\varepsilon}(n))$, then we call $\cla_{-\varepsilon} \in SI(S^{-\varepsilon}(n))$ the {\emph{dual correspondent}} of $\cla_{\varepsilon}$.

\subsection{Preliminaries: the case $\e=1$}\label{Preliminaries}

Let $B$ be a spin block of $S^+(n)$ of weight $w =w(B)>0$ and sign $\delta=\da(B)$, and let $B_0$ be the principal spin block of $S^{\da}(pw)$. Let $r=n-wp$. Let $\mu$ be the $\p$-core of $B$, so that $\sa(\mu)=\da$. The characters in $B$ are indexed by the $\p$-quotients of weight $w$. For any bar-partition $\la$ with $\p$-core $\mu$, we denote the $\p$-quotient of $\la$ by $\la^{(\p)}$ (so that $\sa(\la)=\da \sa(\la^{(\p)})$), and we let $\tilde{\la}$ be the bar-partition of $wp$ with empty $\p$-core, and $\p$-quotient  $\la^{(\p)}$.
\begin{lemma}\label{heightpreserving}If $\da = 1$, then, with the above notation, $\la \longmapsto \tla$ induces a sign-preserving bijection $\cal{I}$ between $B$ and $B_0$ which is also height-preserving. Furthermore,
$$\cyc{\la}(1)_{p'} \equiv \displaystyle \pm \frac{(n!)_{p'}}{((wp)!)_{p'}(r!)_{p'}} \cyc{\mu}(1)_{p'} \cyc{\tla}(1)_{p'} \modp.$$

\end{lemma}

\begin{proof}
Since $\da=1$, we have $\cyc{\tla} \in B_0$ and $\sa(\cyc{\la})=\sa(\la)=\sa(\tla)=\sa(\cyc{\tla})$, so that $\cal{I} \colon \cyc{\la} \longmapsto \cyc{\tla}$ is sign preserving, and therefore gives a bijection between $B$ and $B_0$. 

Now, for any $\cyc{\la} \in B$, we have
$$\cyc{\la}(1)=2^{\floor{(n-m(\la))/2}} \frac{n!}{\bar{h}(\la)} \; \; \mbox{and} \; \; \cyc{\tla}(1)=2^{\floor{(wp-m(\tla))/2}} \frac{(wp)!}{\bar{h}(\tla)}.$$
In particular, since $B$ and $B_0$ have a common defect group $X$ (which is a Sylow $p$-subgroup of $S(pw)$ and $S^+(pw)$), the heights of $\cyc{\la}$ and $\cyc{\tla}$ are
$${\mathfrak{h}}(\cyc{\la})= \nu (|X|) - \nu (\bar{h}(\la)) \; \; \mbox{and} \; \; {\mathfrak{h}}(\cyc{\tla})= \nu (|X|) - \nu (\h(\tla)) \; \; \mbox{respectively}.$$
Now, using the notation of section 2, we have $\h(\la)=\h_{\la,p} \h_{\la,p'}$ and $\h(\tla)=\h_{\tla,p} \h_{\tla,p'}$, so that $\nu (\h(\la))=\nu (\h_{\la,p})$ and $\nu (\h(\tla))=\nu (\h_{\tla,p})$. However, because of the bijection beween bars of length divisible by $p$ in $\la$ and bars in the $\p$-quotient $\la^{(\p)}$, we have $\h_{\la,p}=p^w\h(\la^{(\p)})=p^w\h(\tla^{(\p)})=\h_{\tla,p}$, whence $\nu (\h(\la))=\nu (\h(\tla))$ and ${\mathfrak{h}}(\cyc{\la})={\mathfrak{h}}(\cyc{\tla})$. This proves that $\cal{I}$ is height-preserving. We also get
$$\frac{\cyc{\la}(1)_{p'}}{\cyc{\tla}(1)_{p'}} = \frac{2^{\floor{(n-m(\la))/2}}}{2^{\floor{(wp-m(\tla))/2}}} \frac{(n!)_{p'}}{((wp)!)_{p'}} \frac{(\h_{\tla,p})_{p'} \h_{\tla,p'}}{(\h_{\la,p})_{p'} \h_{\la,p'}}=\frac{2^{\floor{(n-m(\la))/2}}}{2^{\floor{(wp-m(\tla))/2}}} \frac{(n!)_{p'}}{((wp)!)_{p'}} \frac{ \h_{\tla,p'}}{ \h_{\la,p'}}.$$
If we write $\la^{(\p)}=(\la^{(0)}, \, \la^{(1)}, \, \ldots, \, \la^{((p-1)/2)})$, then, by \cite[Corollary (2.6)]{Olsson-FrobeniusSymbols}, we have $m(\la)=m(\la^{(0)})+m(\mu) + 2 a(\la)$ and $m(\tla)=m(\la^{(0)})+m(\emptyset) + 2 a(\tla)$. This implies that $$\floor{(n-m(\la))/2}=\floor{(n-m(\la^{(0)})-m(\mu))/2}-a(\la)$$and$$\floor{(wp-m(\tla))/2}=\floor{(wp-m(\la^{(0)}))/2}-a(\tla).$$Now Proposition \ref{Fong's Lemma} gives $ \h_{\la,p'} \equiv \pm 2^{-a(\lambda)}  \bar{h}(\mu) \modp$ and  $ \h_{\tla,p'} \equiv \pm 2^{-a(\tla)} \modp$. This yields
$$\frac{2^{\floor{(n-m(\la))/2}}}{ \h_{\la,p'}} \equiv \pm \frac{2^{\floor{(n-m(\la^{(0)})-m(\mu))/2}}}{\bar{h}(\mu)} \; \;  \modp
$$and$$ \frac{ \h_{\tla,p'}}{2^{\floor{(n-m(\tla))/2}}} \equiv \pm \frac{1}{2^{\floor{(wp-m(\la^{(0)}))/2}}} \; \;  \modp.$$
By hypothesis, we have $\da=\sa(\mu)=(-1)^{r-m(\mu)}=1$, so that $r-m(\mu)$ is even. Thus
$\floor{(n-m(\la^{(0)})-m(\mu))/2}=\floor{(wp-m(\la^{(0)}))/2}-\floor{(r-m(\mu))/2}$, which in turns implies, together with the above,
$$\frac{\cyc{\la}(1)_{p'}}{\cyc{\tla}(1)_{p'}} \equiv \pm  \frac{(n!)_{p'}}{((wp)!)_{p'}} \frac{2^{\floor{(r-m(\mu))/2}}}{\bar{h}(\mu)}=\pm  \frac{(n!)_{p'}}{((wp)!)_{p'}(r!)_{p'}} \cyc{\mu}(1)_{p'} \; \; \modp. $$

\end{proof}

The corresponding result when $\da=-1$ is given by the following

\begin{lemma}\label{heightreversing}If $\da = -1$, then, with the above notation, $\la \longmapsto \tla$ induces a sign-preserving bijection $\cal{I} \colon \cla \longmapsto \ctla_-$ between $B$ and $B_0$ which is also height-preserving. Furthermore, if $\ctla_+ \in B_0^* \subset SI(S^+(pw))$ is the dual correspondent of $\ctla_-$, then
$$\cyc{\la}(1)_{p'} \equiv \displaystyle \pm \frac{(n!)_{p'}}{((wp)!)_{p'}(r!)_{p'}} \cyc{\mu}(1)_{p'} \cyc{\tla}_+(1)_{p'} 2^{s(\la)} \modp,$$
where $s(\la)=1$ if $\sa(\tla)=-1$ and $s(\la)=0$ if $\sa(\tla)=1$.
\end{lemma}

\begin{proof}
Since $\da=-1$, we have $\sa(\tla)=-\sa(\la)$, and $\sa(\ctla_-)=-\sa(f^-(\ctla_-))=-\sa(\tla)=\sa(\cla)$, so that $\cal{I} \colon \cyc{\la} \longmapsto \cyc{\tla}_-$ is sign preserving, and therefore gives a bijection between $B$ and $B_0$. 

Now, for any $\cyc{\la} \in B$, we have
$$\cyc{\la}(1)=2^{\floor{(n-m(\la))/2}} \frac{n!}{\bar{h}(\la)} \; \; \mbox{and} \; \; \cyc{\tla}_+(1)=2^{\floor{(wp-m(\tla))/2}} \frac{(wp)!}{\bar{h}(\tla)},$$
and, by duality,$$\ctla_-(1)= \left\{ \begin{array}{ll} \ctla_+(1) & \mbox{if} \;  \sa(\ctla_-)=1 \\
\ctla_+(1)/2 & \mbox{if} \;  \sa(\ctla_-)=-1 \end{array} \right. .$$
As in the proof of Lemma \ref{heightpreserving}, this implies that ${\cal I}$ is height-preserving.

As above, we have $$\floor{(n-m(\la))/2}=\floor{(n-m(\la^{(0)})-m(\mu))/2}-a(\la)$$and$$\floor{(wp-m(\tla))/2}=\floor{(wp-m(\la^{(0)}))/2}-a(\tla),$$so that Proposition \ref{Fong's Lemma} yields
$$\frac{2^{\floor{(n-m(\la))/2}}}{ \h_{\la,p'}} \equiv \pm \frac{2^{\floor{(n-m(\la^{(0)})-m(\mu))/2}}}{\bar{h}(\mu)} \; \;  \modp
$$and$$ \frac{ \h_{\tla,p'}}{2^{\floor{(n-m(\tla))/2}}} \equiv \pm \frac{1}{2^{\floor{(wp-m(\la^{(0)}))/2}}} \; \;  \modp.$$
However, this time, we have $\da=\sa(\mu)=(-1)^{r-m(\mu)}=-1$, so that $r-m(\mu)$ is odd. Thus
$$\begin{array}{rcl} \floor{(n-m(\la^{(0)})-m(\mu))/2} & = & \floor{(wp-m(\la^{(0)})+1)/2}-\floor{(r-m(\mu)-1)/2} \\ & = & \floor{(wp-m(\la^{(0)})+1)/2}-\floor{(r-m(\mu))/2}.\end{array}$$Now$$\floor{(wp-m(\la^{(0)})+1)/2}-\floor{(wp-m(\la^{(0)}))/2} = \left\{ \begin{array}{l} 1 \; \mbox{if $wp-m(\la^{(0)})$ is odd}  \\ 0 \; \mbox{if $wp-m(\la^{(0)})$ is even}  \end{array} \right. .$$
But $(-1)^{wp-m(\la^{(0)})}=(-1)^{wp-|\la^{(0)}|}(-1)^{|\la^{(0)}|-m(\la^{(0)})}=(-1)^{wp-|\la^{(0)}|}\sa(\la^{(0)})$, and, since $p$ is odd, $(-1)^{wp-|\la^{(0)}|}=(-1)^{w-|\la^{(0)}|}$, so that $(-1)^{wp-m(\la^{(0)})}=\sa(\la^{(\p)})=-\sa(\la)=\sa(\tla)$. This implies the result.

\end{proof}

We now turn to the $p$-local situation. As mentionned above, the defect group $X$ of the block $B$ can be chosen to be a Sylow $p$-subgroup of $S(pw)$ and of $S^+(pw)$. If we let $N_0(X)=N_{S(pw)}(X)$, then we have$$N_{S(n)}(X)=N_0(X) \times S(r) \; \; \mbox{and} \; \;  N:= N_{S^+(n)}(X)=N_0(X)^+\hat{ \times} S^+(r).$$
In particular, we have, writing ${\cal N}= |N_0(X)^+|=|N_{S(pw)}(X)^+|=|N_{S^+(pw)}(X)|$,
$$[ S^+(n) \colon N_{S^+(n)}(X) ] = \frac{ 2 n!}{ | N_0(X)^+ \hat{\times} S^+(r) | } =\frac{2n!}{({\cal N}2r!)/2}  =\frac{2n!}{{\cal N}r!}$$and$$[ S^+(pw) \colon N_{S^+(pw)}(X) ] =\frac{2(pw)!}{{\cal N}}.$$
The following result will also be useful later; note that, if $X$ is a Sylow $p$-subgroup of $S^+(pw)$, then it is also a Sylow $p$-subgroup of $S^-(pw)$ (since $S^-(pw)$ is of index 2 in $S^+(pw)$, while $p$ is odd).

\begin{lemma}\label{normalizer}
If $X$ is a Sylow $p$-subgroup of $S^+(pw)$, then$$[S^-(pw) \colon N_{S^-(pw)}(X)]=[S^+(pw) \colon N_{S^+(pw)}(X)] $$and$$ [S^-(n) \colon N_{S^-(n)}(X)]=[S^+(n) \colon N_{S^+(n)}(X)].$$
\end{lemma}

\begin{proof}

We start by noticing that we must have $$  \frac{[S^+(pw) \colon N_{S^+(pw)}(X)]}{[S^-(pw) \colon N_{S^-(pw)}(X)] } =  \frac{2}{ [N_{S^+(pw)}(X) \colon N_{S^-(pw)}(X)]} \in \{1, \, 2 \}$$
(since $[N_{S^+(pw)}(X) \colon N_{S^-(pw)}(X)]\in \{1, \, 2 \}$).

However, by Sylow's Theorems, we have $[S^-(pw) \colon N_{S^-(pw)}(X)]  \equiv 1 \modp$. Thus $[S^+(pw) \colon N_{S^+(pw)}(X)] =2 [S^-(pw) \colon N_{S^-(pw)}(X)] $ would yield \\$[S^+(pw) \colon N_{S^+(pw)}(X)]  \equiv 2 \modp$, a contradiction, $X$ being a Sylow $p$-subgroup of $S^+(pw)$. Hence $[S^-(pw) \colon N_{S^-(pw)}(X)]=[S^+(pw) \colon N_{S^+(pw)}(X)] $.

\smallskip
Now $N_{S^+(n)}(X)=N_{S^+(pw)}(X) \hat{\times} S^+(r)$. Hence
$$\begin{array}{rcl} [S^+(n) \colon N_{S^+(n)}(X)] & = & [S^+(n) \colon S^+(pw) \hat{\times} S^+(r)] \, [S^+(pw) \hat{\times} S^+(r) \colon N_{S^+(n)}(X)] \\ & = & [S^+(n) \colon S^+(pw) \hat{\times} S^+(r)] \, [S^+(pw) \colon N_{S^+(pw)}(X)] \\ & = &\displaystyle \frac{2 n!}{ (2 (pw)! 2 r!)/2} [S^+(pw) \colon N_{S^+(pw)}(X)]  \\ &  = & \displaystyle \frac{n!}{  (pw)!  r!} [S^-(pw) \colon N_{S^-(pw)}(X)] . \end{array}$$
On the other hand, we have, {\bf{as sets}}, a disjoint union
$$N_{S^-(n)}(X)=(N_{S^-(pw)}(X) \hat{\times} S^-(r) )   \bigcup   (  N_{S^+(pw) \setminus S^-(pw)}(X) \hat{\times} (S^+(r) \setminus S^-(r) ) ),$$
so that $|N_{S^-(n)}(X)|=2 |N_{S^-(pw)}(X) \hat{\times} S^-(r) |$. We thus get
$$\begin{array}{rcl} [S^-(n) \colon N_{S^-(n)}(X)] & = & \displaystyle \frac{|S^+(n)|}{|S^-(pw) \hat{\times} S^-(r)|} \frac{|  S^-(pw) \hat{\times} S^-(r) |}{2 | N_{S^-(pw)}(X) \hat{\times} S^-(r)|}  \\ & = & \displaystyle\frac{ n!}{ ((pw)! r!)/2} \frac{ [S^-(pw) \colon N_{S^-(pw)}(X)]}{2}  \\ &  = &\displaystyle \frac{n!}{  (pw)!  r!} [S^-(pw) \colon N_{S^-(pw)}(X)], \end{array}$$
whence the result.

\end{proof}

If $\mu \succ r $ is the $\p$-core of $B$, we choose $\ga \in SI(S^+(r))$ such that $f^+(\ga)=\mu$. We have $\ga=\ga^a$ if and only if $\da=1$.

If we then denote by $b$ the Brauer correspondent of $B$ in $N$, we have (see the proof of \cite[Theorem 2.2]{Michler-Olsson})
$$b= \{ \chi \hat{\otimes} \ga, \,  \chi \hat{\otimes} \ga^a \; | \; \chi \in SI(N_0(X)^+) \}= \{ \chi \hat{\otimes} \ga, \,  \chi \hat{\otimes} \ga^a \; | \; \chi \in \ba_0 \},$$
where $\ba_0$ is {\bf{the}} spin block of $N_0(X)^+$, and thus the Brauer correspondent of the principal spin block of $S^+(pw)$. In particular, $\ba_0=b_0$ if $\da=1$ and $\ba_0=b_0^*$ if $\da=-1$.

For any $\chi \in \ba_0$, we have $ (\chi \hat{\otimes} \ga) (1) = ( \chi \hat{\otimes} \ga^a ) (1)= 2^{\floor{s/2}} \chi (1) \ga (1)$, where $s$ is the number of n.s.a. characters in $\{ \chi , \, \ga \}$. If $\da=1$, we therefore get $s=0$ (if $\chi^a=\chi$) or $s=1$ (if $\chi^a \neq \chi$), so that $\floor{s/2}=0$ and  $ (\chi \hat{\otimes} \ga) (1) =  \chi (1)  \ga (1)$. If $\da = -1$, we have $s=1$ and $\floor{s/2}=0$ if $\chi^a=\chi$, and $s=2$ and $\floor{s/2}=1$ if $\chi^a \neq \chi$, so that
$$ (\chi \hat{\otimes} \ga) (1) = ( \chi \hat{\otimes} \ga^a ) (1)= \left\{ \begin{array}{l}  \chi (1)  \ga (1) \;  \mbox{if} \;  \chi^a = \chi \\  2 \chi (1)  \ga (1) \;  \mbox{if} \;  \chi^a \neq \chi \end{array} \right. .$$

\subsection{Reduction Theorem}
We can now prove the main result of this section:

\begin{theorem}\label{Reduction Theorem}
Let $B$ be a spin block of $S^{\e}(n)$ of weight $w= w(B) >0$ and sign $\da=\da(B)$ and let $b$ be its Brauer correspondent in $N_{S^{\e}(n)}(X)$, where $X$ is a defect group of $B$. Suppose the Isaacs-Navarro Conjecture holds for the principal spin block of $S^+(pw)$ via a sign-preserving bijection. Then it also holds for $B$.
\end{theorem} 

\begin{proof}
We first suppose ${\mathbf{\e=1}}$. 

\noindent
We use the same notation as in Section \ref{Preliminaries}. Let $B_0$ be the principal spin block of $S^{\da}(pw)$ and $b_0$ be its Brauer correspondent. Let $\mu \succ r=n-wp$ be the $\p$-core of $B$, and $\ga \in SI(S^+(r))$ such that $f^+(\ga)=\mu$. If $\la$ is a bar-partition of $n$ with $\p$-core $\mu$ and $\p$-quotient $\la^{(\p)}$, let $\tilde{\la}$ be the bar-partition of $wp$ with empty $\p$-core and $\p$-quotient  $\la^{(\p)}$.

Suppose furthermore that $\da = 1$. Then, by Lemma \ref{heightpreserving}, $\la \longmapsto \tla$ induces a sign-preserving bijection $\cal{I}$ between $B$ and $B_0$ which is also height-preserving, and
$$\cyc{\la}(1)_{p'} \equiv \displaystyle \pm \frac{(n!)_{p'}}{((wp)!)_{p'}(r!)_{p'}} \ga(1)_{p'} \cyc{\tla}(1)_{p'} \modp.$$
Now let $c=[S^+(pw) \colon N_{S^+(pw)}(X) ]_{p'}$, and let $\varphi \colon M(B_0) \longrightarrow M(b_0)$ be a bijection such that, for each $k$ such that $(p, \, k)=1$, we have $M_k(b_0)= \varphi(M_{ck}(B_0))$ (such a $\varphi$ exists by hypothesis). We have $b = \{ \chi \hat{\otimes} \ga \; | \; \chi \in b_0 \}$, and, by \cite[Proposition 1.2]{Michler-Olsson}, $ \chi \hat{\otimes} \ga =  \psi \hat{\otimes} \ga $ if and only if $\psi \in \{ \chi, \, \chi^a \}$ and ( $\sa (\chi) \sa(\ga) = \sa (\chi)=1$) or ($ \sa ( \chi) =-1$ and $ \chi = \psi $), i.e. $ \chi \hat{\otimes} \ga =  \psi \hat{\otimes} \ga $ if and only if $\psi = \chi$. Thus, by the results of Section \ref{Preliminaries}, $\chi \longmapsto \chi \hat{\otimes} \ga$ is a height-preserving bijection between $b_0$ and $b$. Hence
$$\Phi \colon \left\{ \begin{array}{rcl} M(B) & \longrightarrow & M(b) \\ \cyc{\la} & \longmapsto & \varphi(\cyc{\tla}) \hat{\otimes} \ga \end{array} \right.$$is a (height-preserving) sign-preserving bijection.

Now, if $\cla \in M(B)$, then $\Phi(\cla)(1)=(\varphi(\ctla))(1)\ga(1)$, so that
$$ \begin{array}{rcl} \Phi(\cla)(1)_{p'} & = & (\varphi(\ctla))(1)_{p'}\ga(1)_{p'} \\ & \equiv &\displaystyle \pm \frac{\ctla (1)_{p'}}{c} \ga(1)_{p'} \modp \; \; (\mbox{by definition of} \; \varphi) \\  & \equiv & \displaystyle \pm \frac{1}{c} \left( \frac{r!(pw)!}{n!} \right)_{p'} \cla (1)_{p'} \modp, \end{array}$$
and
$$\begin{array}{rcl} \displaystyle \frac{1}{c} \left( \frac{r!(pw)!}{n!} \right)_{p'} & = & \displaystyle \frac{1}{[S^+(pw) \colon N_{S^+(pw)}(X) ]_{p'}} \left( \frac{r!(pw)!}{n!} \right)_{p'} \\ & = &  \displaystyle \left( \frac{\cal N}{2(pw)!} \right)_{p'} \left( \frac{r!(pw)!}{n!} \right)_{p'}=\left( \frac{{\cal N}r!}{2n!} \right)_{p'} \\ & = & \displaystyle \frac{1}{[S^+(n) \colon N_{S^+(n)}(X) ]_{p'}}, \end{array}$$
whence we finally get
$$ \Phi(\cla)(1)_{p'} \equiv  \pm  \displaystyle \frac{\cla (1)_{p'}}{[S^+(n) \colon N_{S^+(n)}(X) ]_{p'}} \modp,$$i.e. $\Phi$ is an Isaacs-Navarro bijection between $B$ and $b$.

\medskip

Suppose now that $\da=-1$. Then $B_0$ is the principal spin block of $S^-(pw)$, its dual $B_0^*$ is the principal spin block of $S^+(pw)$, and $b_0^*$ is the Brauer correspondent of $B_0^*$. Writing $D_+$ for the set of s.a. characters in $D$ and $D_-$ for the set of pairs of n.s.a. characters in $D$ (so that $|M_k(B_0)|=|M_k(B_0)_+|+2|M_k(B_0)_-|$), we thus have, for each $k$ such that $(p, \, k)=1$, the following equalities:
$$|M_k(B_0)_+|=|M_k(B_0^*)_-|=|M_{k/c}(b_0^*)_-| $$and$$ |M_k(B_0)_-|=|M_{2k}(B_0^*)_+|=|M_{2k/c}(b_0^*)_+|,$$where $c=[S^+(pw) \colon N_{S^+(pw)}(X) ]_{p'}=[S^-(pw) \colon N_{S^-(pw)}(X) ]_{p'}$ (by Lemma \ref{normalizer}).

On the other hand, we have $b=\{ \chi \hat{\otimes} \ga, \, \chi \hat{\otimes} \ga^a \, | \, \chi \in \ba_0=b_0^* \}$. For any $\chi, \, \psi \in \ba_0$ and $\ga_1, \, \ga_2 \in \{ \ga, \, \ga^a \}$, we have $\chi  \hat{\otimes} \ga_1 = \psi \hat{\otimes} \ga_2$ if and only if $\chi$ and $\psi$ are associate and
$$\left\{ \begin{array}{l} \sa(\chi) \sa(\ga_1)= - \sa(\chi) = 1  \\ \mbox{or} \\ \sa(\chi) \sa(\ga_1)= - \sa(\chi) = -1 \; \mbox{and} \; [ (\chi = \psi, \, \ga_1=\ga_2) \; \mbox{or} \;  (\chi \neq \psi, \, \ga_1 \neq \ga_2) ]\end{array} \right. $$
Hence, if $\chi \in \ba_{0+}$, then we get two irreducible characters, $\chi \hat{\otimes} \ga$ and $\chi \hat{\otimes} \ga^a$, while, if $\chi \in \ba_{0-}$, then we get one irreducible character, $\chi \hat{\otimes} \ga = \chi^a \hat{\otimes} \ga=\chi^a \hat{\otimes} \ga^a = \chi \hat{\otimes} \ga^a$. Note that $\chi \longmapsto \chi \hat{\otimes} \ga$ and $\chi \longmapsto \chi \hat{\otimes} \ga^a$ are height preserving. Using the equalities above, as well as Lemma \ref{heightreversing}, we obtain the following height-preserving and sign-preserving bijection:
$$\Phi \colon \left\{ \begin{array}{ccccccccc} M(B) & {\displaystyle{\mathop{\rightarrow }^+}} &  M(B_0) &  {\displaystyle{\mathop{\rightarrow }^-}}  & M(B_0^*) & {\displaystyle{\mathop{\rightarrow }^+}}  & M(b_0^*) &   {\displaystyle{\mathop{\rightarrow }^-}} & M(b) \\
\cla & \mapsto & \ctla_- & \mapsto & \ctla_+ & \mapsto & \varphi (\ctla_+) & \mapsto & \varphi (\ctla_+) \hat{\otimes} \cyc{\mu}
\end{array} \right. ,$$
where, as before, $\varphi$ is the (sign-preserving) Isaacs-Navarro bijection we supposed exists between $ M(B_0^*) $ and $ M(b_0^*) $, and $ {\displaystyle{\mathop{\rightarrow }^+}}$ (respectively  ${\displaystyle{\mathop{\rightarrow }^-}}$) denotes a sign-preserving (respectively sign-inversing) bijection.

Now, by hypothesis, $\ctla_+ (1)_{p'} \equiv c \varphi(\ctla_+) (1)_{p'} \modp$, so that, by Lemma \ref{heightreversing}, we obtain
$$\cla (1)_{p'} \equiv \left( \frac{n!}{(wp)!r!} \right)_{p'} c \varphi(\ctla_+) (1)_{p'} \cyc{\mu} (1)_{p'} 2^{s(\la)} \modp ,$$
and, as in the case $\da=1$, we have $\left( \frac{n!}{(wp)!r!} \right)_{p'} c = [ S^+(n) \colon N_{S^+(n)}(X) ]_{p'}$. Finally, since $\sa(\tla)=1 \Longleftrightarrow \sa(\ctla_+)=1 \Longleftrightarrow \ctla_+ \; \mbox{is s.a.} \;  \Longleftrightarrow \varphi (\ctla_+) \; \mbox{is s.a.}$, we get
$$  \begin{array}{rcl} \varphi(\ctla_+) (1)_{p'} \cyc{\mu} (1)_{p'} 2^{s(\la)} & = & \left\{ \begin{array}{ll}  2 \varphi(\ctla_+) (1)_{p'} \cyc{\mu} (1)_{p'} & \mbox{if} \; \sa(\tla)=-1 \\ \varphi(\ctla_+) (1)_{p'} \cyc{\mu} (1)_{p'} & \mbox{if} \; \sa(\tla)=1 \end{array} \right. \\ & = & ( \varphi(\ctla_+) \hat{\otimes} \cyc{\mu}) (1)_{p'} \\ & = & \Phi (\cla)(1)_{p'}, \end{array} $$
whence $\cla (1)_{p'} \equiv [ S^+(n) \colon N_{S^+(n)}(X) ]_{p'} \Phi (\cla)(1)_{p'} \modp$, i.e. $\Phi$ is a (sign-preserving) Isaacs-Navarro bijection between $M(B)$ and $M(b)$.

\medskip
\noindent
We now suppose ${\mathbf{\e=-1}}$.

\noindent
In this case, $B$ is a spin block of $S^-(n)$ and $b$ is its Brauer correspondent in $N_{S^-(n)}(X)$. Thus $B^*$ is a spin block of $S^+(n)$, and, by \cite[Lemma 2.3]{Michler-Olsson} (which is due to H. Blau), the dual $b^*$ of $b$ is the Brauer correspondent of $B^*$. By the case $\e=1$, there exists a sign-preserving Isaacs-Navarro bijection $\varphi \colon M(B^*) \longrightarrow M(b^*)$. We define the sign-preserving bijection$$\Phi \colon \left\{ \begin{array}{rcl} M(B) & \longrightarrow & M(b) \\ \cla & \longmapsto & (\varphi(\cla^*))^* \end{array} \right.$$
By Lemma \ref{normalizer}, we have $c = [S^-(n) \colon N_{S^-(n)}(X)]_{p'}=[S^+(n) \colon N_{S^+(n)}(X)]_{p'}$, and, for each $k$ such that $(p, \, k)=1$, we have
$$|M_{ck}(B)_+|=|M_{ck}(B^*)_-|=|M_{k}(b^*)_-|=|M_k(b)_+| $$and$$|M_{ck}(B)_-|=|M_{2ck}(B^*)_+|=|M_{2k}(b^*)_+|=|M_k(b)_-|,$$
whence $|M_{ck}(B)|=|M_{ck}(B)_+| + 2 |M_{ck}(B)_-|=|M_{k}(b)|$.

\end{proof}

\section{Principal block}

By Theorem \ref{Reduction Theorem}, it is now sufficent to prove that the Isaacs-Navarro Conjecture holds for the principal spin block of $S^+(pw)$ via a sign-preserving bijection. Throughout this section, we therefore consider the following situation. We take $G=S^+(pw)$ (where $w\geq 1$ is an integer), $B$ the principal spin block of $G$, and $b$ the Brauer correspondent of $B$. Hence $b$ is the principal spin block of $N_G(X)$ for some $X \in Syl_p(G)$.

\subsection{Spin characters of height 0 of the normalizer}

The normalizer $N^+=N_G(X)$ and its irreducible spin characters are described in Sections 3 and 4 of \cite{Michler-Olsson}. Let $pw=\sum_{i=1}^k t_i p^i$ be the $p$-adic decomposition of $pw$. We then have $N^+=[N_1 \wr S(t_1)]^+ \hat{\times} \cdots \hat{\times} [N_k \wr S(t_k)]^+$, where, for each $1 \leq i \leq k$, $N_{i}=N_{S(p^{i})}(X_{i})$ for some $X_{i} \in Syl_p(S(p^{i}))$.

Now fix $1 \leq i \leq k$, and let $e_{i}=(p^{i}-1)/2$. Then $H_i^+=(N_{i} \wr S(t_{i}))^+=M_i^+ S_{t_i}^+$, where $M_i^+=N_{i}^{(1)+} \hat{\times} \cdots \hat{\times} N_{i}^{(t_i)+} \triangleleft   H_i^+$ and $S_{t_i} \cong \Delta_{p^{i}}S(t_i) \subset S(p^{i}t_i)$, and where
$$S_{t_i}^+ \cong \left\{ \begin{array}{ll} \hat{S}(t_i) & \mbox{if} \; p^{i} \equiv 1 \modquat \\  \tilde{S}(t_i) &\mbox{if} \; p^{i} \equiv -1 \modquat \end{array} \right. .$$
By \cite[Proposition 3.9]{Michler-Olsson}, $N_{i}^+$ has one s.a. spin character $\za_0$ of degree $(p-1)^{i}$, and $e_i=(p^{i}-1)/2$ pairs of n.s.a. spin characters $\{ \za_1, \, \za_1^a, \, \ldots , \, \za_{e_i}, \, \za_{e_i}^a \}$ of degree 1.

Let ${\cal A}_{i} = \{ (t_i^{(0)},  \, t_i^{(1)}, \, \ldots , \, t_i^{(e_i)} ) \, | \, t_i^{(j)} \in \N \cup \{0\}, \, \sum_{j=0}^{e_i} t_i^{(j)} = t_{i} \}$.  Then, by \cite[Proposition 3.12]{Michler-Olsson}, a complete set of representatives for the $S_{t_i}^+$-conjugacy classes in $SI_0(M_i^+)$ is given by
$${\cal R} = \{ \ta_{\bf s} \, | \, {\bf s} \in {\cal A}_{i} \} \cup \{ \ta_{\bf s}^a \, | \, {\bf s}= (t_i^{(0)},  \, t_i^{(1)}, \, \ldots , \, t_i^{(e_i)} )  \in {\cal A}_{i}, \, t_i-t_i^{(0)} \, \mbox{odd},  \, t_i^{(0)} \leq 1 \},$$
where $\ta_{\bf s} = \ta_0 \hat{\otimes} \ta_1 \hat{\otimes} \cdots \hat{\otimes} \ta_{e_i}$, with $\ta_j=\za_j  \hat{\otimes} \cdots \hat{\otimes} \za_j$ ($t_i^{(j)}$ factors). Also, the inertial subgroup $T_i^+=I_{H_i^+}(\ta_{\bf s})$ of $\ta_{\bf s}$ in $H_i^+ $ satisfies
$$T_i^+ / M_i^+ \cong \left\{ \begin{array}{ll} A(t_i^{(0)}) \times S(t_i^{(1)}) \times \cdots \times S(t_i^{(e_i)}) & \mbox{if} \; t_i-t_i^{(0)} \; \mbox{is odd} \\ S(t_i^{(0)}) \times S(t_i^{(1)}) \times \cdots \times S(t_i^{(e_i)}) & \mbox{if} \; t_i-t_i^{(0)} \; \mbox{is even} \end{array} \right. .$$ 
We can now describe how to induce each $\ta_j$ from $M_i^{(j)+}=(N_{i}^+)^{\hat{\times}t_i^{(j)}}$ to the corresponding factor $T_i^{(j)+}$ of its inertial subgroup.

\begin{proposition}\cite[Proposition 4.4]{Michler-Olsson}\label{Prop 4.4}
If $\za$ is a n.s.a. linear representation of $N_{i}^+$, then $\ta_j=\za^{t_i^{(j)}}= \za \hat{\otimes} \cdots \hat{\otimes} \za \in \Irr (M_i^{(j)+})$ can be extended to a negative representation $D_{\za} \in \Irr (T_i^{(j)+})$, and every irreducible constituent $V$ of $\ta \uparrow^{T_i^{(j)+}}$ is of the form $V = D_{\za} \otimes R$, where $R$ is an irreducible representation of $T_i^{(j)+}/M_i^{(j)+} \cong S(t_i^{(j)})$. If $t_i^{(j)}$ is odd, then every irreducible constituent $V$ of $\ta \uparrow^{T_i^{(j)+}}$ is n.s.a., and, if $t_i^{(j)}$ is even, then every irreducible constituent $V$ of $\ta \uparrow^{T_i^{(j)+}}$ is s.a..

\end{proposition}

In the above notation, if $\psi$ is the character of $V = D_{\za} \otimes R$, and if $R$ has character $\chi_{\la} \in \Irr(S(t_i^{(j)}))$, then $\psi(1)=\za^{t_i^{(j)}}(1) \chi_{\la}(1)$. Also, since $\za$ is n.s.a., we have $\za^{t_i^{(j)}}(1) = 2^{\floor{t_i^{(j)}/2}} \za(1)^{t_i^{(j)}}=2^{\floor{t_i^{(j)}/2}}$, and $\psi(1)=2^{\floor{t_i^{(j)}/2}} \chi_{\la}(1)$. Finally, $\psi$ is s.a. if and only if $t_i^{(j)}$ is even.

\begin{proposition}\cite[Proposition 4.8]{Michler-Olsson}\label{Prop 4.8}
Let $t_i^{(0)} \geq 4$, and let $D$ be the s.a. spin representation of $N_{i}^+$ with degree $(p-1)^{i}$.
Then $D^{t_i^{(0)}}= D \hat{\otimes} \cdots \hat{\otimes} D \in \Irr (M_i^{(0)+})$ can neither be extended to an irreducible representation of $T_i^{(0)-}=M_i^{(0)+}A_{t_i^{(0)}}^+$ nor to one of $T_i^{(0)+}=M_i^{(0)+}S_{t_i^{(0)}}^+$. Furthermore, every irreducible constituent $V$ of $D^{t_i^{(0)}} \uparrow^{T_i^{(0)+}}$ is of the form $V = D^{t_i^{(0)}} \otimes S$, where $S$ is an irreducible spin representation of $S^+_{t_i^{(0)}}$, and every irreducible constituent $V$ of $D^{t_i^{(0)}} \uparrow^{T_i^{(0)-}}$ is of the form $V = D^{t_i^{(0)}} \otimes S$, where $S$ is an irreducible spin representation of $A^+_{t_i^{(0)}}$. In each case, $V$ is s.a. if and only if $S$ is s.a..
\end{proposition}
In this notation, if $\psi$ is the character of $V $, and if $S$ has character $\chi_S $, then $\psi(1)=\za_0^{t_i^{(0)}}(1) \chi_S(1)$. And, since $\za_0$ is s.a., we have $\psi(1)=(p-1)^{i t_i^{(0)}} \chi_S(1)$.

\medskip
We can now describe all the characters of height 0 in $b$. Recall that these are exactly the spin characters with $p'$-degree in $N^+$. Still writing $pw=\sum_{i=1}^k t_i p^i$ the $p$-adic decomposition of $pw$, Olsson proved in \cite{olsson-blocks} that, for any sign $\sa$,
$$|M(b)_{\sa}| = \displaystyle \sum_{ \{ (\sa_1, \, \ldots , \, \sa_k) \} } \prod_{i=1}^k q^{\sa_i} (\p^i, t_i),$$
where $(\sa_1, \, \ldots , \, \sa_k)$ runs through all $k$-tuples of signs satisfying $\sa_1 \ldots \sa_k=\sa$, and where $q^{\sa_i} (\p^i, t_i)$ denotes the number of all $\p^i$-quotients with sign $\sa_i$ and weight $t_i$.

The correspondence goes as follows. For each $1 \leq i \leq k$, pick ${\bf s_i} \in {\cal A}_i$ and the corresponding $\ta_{\bf s_i} \in SI_0(M_{\bf s_i}^+)$ (where, if ${\bf s_i}= (t_i^{(0)}, \, t_i^{(1)}, \, \ldots , \, t_i^{(e_i)})$, then $M_{\bf s_i}^+ = (N_i^{\hat{\times} t_i^{(0)}})^+ \hat{\times} \cdots \hat{\times} (N_i^{\hat{\times} t_i^{(e_i)}})^+$). Inducing $\ta_{\bf{ s_i}}$ (or $\ta_{\bf{ s_i}} + \ta_{\bf{ s_i}}^a$ if $t_i - t_i^{(0)}$ is odd and $t_i^{(0)} \leq 1$) to its inertial subgroup $T_i^+$, we obtain s.a. irreducible constituents and pairs of n.s.a. irreducible constituents described by Propositions \ref{Prop 4.4} and \ref{Prop 4.8} and labeled by the $\p^i$-quotients of weight $t_i$: if $Q_i=(\la_i^{(0)}, \, \la_i^{(1)}, \, \ldots , \, \la_i^{(e_i)})$ is a $\p^i$-quotient of weight $t_i$, then $\cyc{\Psi_{Q_i}} = \cyc{\psi_i^{(0)}} \hat{\times} \cyc{\psi_i^{(1)}} \hat{\times} \cdots \hat{\times} \cyc{\psi_i^{(e_i)}} \in SI_0(T_i^+)$. Also, by Propositions \ref{Prop 4.4} and \ref{Prop 4.8}, 

\begin{itemize}
\item{} For $1 \leq j \leq e_i$, $\psi_i^{(j)} (1)=2^{\floor{t_i^{(j)}/2}} \chi_{\la_i^{(j)}}(1)$ (with $\chi_{\la_i^{(j)}} \in \Irr(S(t_i^{(j)}))$) and $\psi_i^{(j)} $ is s.a. if and only if $t_i^{(j)}$ is even.
\item{} $\psi_i^{(0)} (1)=(p-1)^{it_i^{(0)}} \chi_{\la_i^{(0)}}(1)$ (with $\chi_{\la_i^{(0)}} \in SI(S^+_{t_i^{(0)}})$ if $t_i-t_i^{(0)}$ is even and $\chi_{\la_i^{(0)}} \in SI(A^+_{t_i^{(0)}})$ if $t_i-t_i^{(0)}$ is odd) and $\psi_i^{(0)} $ is s.a. if and only if $\chi_{\la_i^{(0)}}$ is s.a..
\end{itemize}

Finally,  $\cyc{\Psi_{Q_i}}(1) = 2^{\floor{S_i/2}} \cyc{\psi_i^{(0)}}(1)  \cyc{\psi_i^{(1)}}(1)  \ldots  \cyc{\psi_i^{(e_i)}}(1)$, where $S_i$ is the number of (pairs of) n.s.a. characters in $\{\cyc{\psi_i^{(0)}}, \,  \cyc{\psi_i^{(1)}}, \,   \ldots , \,  \cyc{\psi_i^{(e_i)}}\}$.

\medskip

Inducing to $H_i^+=[N_i \wr S(t_i)]^+$, we obtain a s.a. irreducible spin character, or a pair of associate (n.s.a.) spin characters, $\cQi$, labeled by $Q_i$.

Given the structure of $T_i^+$, we see that $\sa(\cQi)=(-1)^{t_i-t_i^{(0)}} \sa(\Psi_{Q_i} )$. However, we have $\sa({\Psi_{Q_i}} )=  \sa(\psi_i^{(0)}) \sa(\psi_i^{(1)}) \ldots \sa(\psi_i^{(e_i)})$. Also, for $1 \leq j \leq e_i$, we have $\sa(\psi_i^{(j)})=(-1)^{t_i^{(j)}}$, and $\sa(\psi_i^{(0)})=\sa(\chi_{\la_i^{(0)}})=\sa(\la_i^{(0)}) (-1)^{t_i-t_i^{(0)}}$, so that $\sa({\Psi_{Q_i}} )= \sa(\la_i^{(0)})$ (since $\sum_{j=1}^{e_i} t_i^{(j)}=t_i-t_i^{(0)}$) and $\sa(\cQi)=(-1)^{t_i-t_i^{(0)}} \sa(\la_i^{(0)}) = \sa(Q_i)$. Note that, writing $m_i^{(0)}$ for the number of (non-zero) parts in $\la_i^{(0)}$, we have $\sa(\la_i^{(0)}) = (-1)^{t_i^{(0)}-m_i^{(0)}}$, so that  $\sa(\cQi)=(-1)^{t_i-m_i^{(0)}}$, and $\cQi$ is s.a. if and only if ${t_i-m_i^{(0)}}$ is even.

Also, we have $\cQi(1)=(|H_i^+|/|T_i^+|) \cyc{\Psi_{Q_i}}(1)$, unless $\chi_{\la_i^{(0)}}$ is a s.a. irreducible spin character of $A^+_{t_i^{(0)}}$ (i.e. $t_i-t_i^{(0)}$ is odd and $\chi_{\la_i^{(0)}}$ is s.a.), in which case $\cQi(1)=(|H_i^+|/|T_i^+|) \cyc{\Psi_{Q_i}}(1)/2$.

\medskip
Finally, the irreducible characters of height 0 in $b$ are parametrized by the sequences $(Q_1, \, \ldots , \, Q_k)$, where $Q_i$ is a $\p^i$-quotient of weight $t_i$. We have \\$\cyc{ (Q_1, \, \ldots , \, Q_k) }= \cyc{Q_1} \hat{\otimes} \cdots \hat{\otimes} \cyc{ Q_k}$, and $\cyc{ (Q_1, \, \ldots , \, Q_k) }(1)= 2^{\floor{s/2}} \prod_{i=1}^k \cyc{Q_i}(1)$, where $s$ is the number of (pairs of) n.s.a. characters in $\{ \cyc{Q_1}, \, \ldots , \, \cyc{Q_k} \}$. By the above remark on the sign of $\cQi$, we see that $s=| \{ 1 \leq i \leq k \, ; \;  t_i - m_i^{(0)} \; \mbox{odd} \} |$.

\begin{proposition}\label{degspinb}
With the above notation, we have
$$\frac{|N_G(X)|_{p'}}{\cyc{ (Q_1, \, \ldots , \, Q_k) }(1)_{p'}} \equiv \pm \frac{2}{2^{\floor{s/2}}} \prod_{i=1}^k \frac{1}{2^{\floor{ (t_i - m_i^{(0)})/2 }}} \h(Q_i) \modp ,$$
where $s=| \{ 1 \leq i \leq k \, ; \;  t_i - m_i^{(0)} \, \mbox{odd} \} |$, and, for each $1 \leq i \leq k$, $\h(Q_i)$ is the product of all bar-lengths in $Q_i$.

\end{proposition}

\begin{proof}

We have 
$$\frac{|N_G(X)|_{p'}}{\cyc{ (Q_1, \, \ldots , \, Q_k) }(1)_{p'}} = \frac{ \prod_{i=1}^k |H_i^+|}{2^{k-1} \cyc{ (Q_1, \, \ldots , \, Q_k) }(1)_{p'}} = \frac{ \prod_{i=1}^k |H_i^+|}{2^{k-1}  2^{\floor{s/2}} \prod_{i=1}^k \cyc{Q_i}(1)_{p'}}.$$
This gives
$$\frac{|N_G(X)|_{p'}}{\cyc{ (Q_1, \, \ldots , \, Q_k) }(1)_{p'}} =  \frac{1}{2^{k-1} } \frac{1}{2^{\floor{s/2}} } D_1^{(+1)} D_0^{(+1)} D_1^{(-1)} D_0^{(-1)},$$
where, for $\varepsilon \in \{+1 , \,  -1\}$ and $a \in \{0, \, 1 \}$,
$$D_a^{(\varepsilon)} = \prod_{ {\tiny{ \begin{array}{c} 1 \leq i \leq k \\  t_i-t_i^{(0)} \equiv a \modtwo   \\  \sa(\chi_{\la_i^{(0)}})=  \varepsilon \end{array}}} } \frac{|H_i^+|}{\cQi (1)_{p'} }.
$$
Now we have, whenever $(\varepsilon, \, a) \in \{ (+1, \, 0), \, (-1, \, 1), \, (-1, \, 0) \}$,
$$D_a^{(\varepsilon)} = \prod_{ {\tiny{ \begin{array}{c} 1 \leq i \leq k \\  t_i-t_i^{(0)} \equiv a \modtwo   \\  \sa(\chi_{\la_i^{(0)}})=  \varepsilon \end{array}}} }  \frac{|T_i^+|}{\cyc{\Psi_{Q_i}} (1)_{p'} }=
 \prod_{ {\tiny{ \begin{array}{c} 1 \leq i \leq k \\  t_i-t_i^{(0)} \equiv a \modtwo   \\  \sa(\chi_{\la_i^{(0)}})=  \varepsilon \end{array}}} }  \frac{|T_i^+|}{2^{\floor{S_i/2}} \prod_{j=0}^{e_i} \psi_i^{(j)}(1)},$$
while
$$D_1^{(+1)} = \prod_{ {\tiny{ \begin{array}{c} 1 \leq i \leq k \\  t_i-t_i^{(0)} \, \mbox{odd}   \\  \sa(\chi_{\la_i^{(0)}})=  1 \end{array}}} }  \frac{|T_i^+|}{\cyc{\Psi_{Q_i}} (1)_{p'} / 2 }=
 \prod_{ {\tiny{ \begin{array}{c} 1 \leq i \leq k \\  t_i-t_i^{(0)} \, \mbox{odd}  \\  \sa(\chi_{\la_i^{(0)}})=  1 \end{array}}} }  \frac{|T_i^+|}{2^{\floor{S_i/2}-1} \prod_{j=0}^{e_i} \psi_i^{(j)}(1)},$$
where $S_i$ is the number of (pairs of) n.s.a. characters in $\{ \cyc{\psi_i^{(0)}}, \,   \ldots , \,  \cyc{\psi_i^{(e_i)}} \}$.

\medskip
For each $1 \leq i \leq k$, we have $|T_i^+|_{p'}= | T_i^+ / M_i^+ |_{p'} |M_i^+|_{p'}$. Also, $M_i^+ \cong (N_i^+) ^{\hat{\times} t_i}$, so that $|M_i^+|= \frac{|N_i^+|^{t_i}}{2^{t_i-1}}=2 |N_i|^{t_i}$. But $N_i=N_{S(p^i)}(X_i)$ for some $X_i \in Syl_{p}(S(p^i))$; thus $|N_i|=|X_i|. |N_i / X_i|$, and we have $N_i/X_i=K_i \cong (\Z / (p-1)\Z)^i $ (see \cite[page 89]{Michler-Olsson}). Hence $|N_i|_{p'}=(p-1)^i$, $|M_i^+|_{p'}=2(p-1)^{it_i}$, and
$$|T_i^+|_{p'} \equiv \left\{ \begin{array}{cl} (-1)^{i t_i} t_i^{(0)} ! t_i^{(1)} ! \ldots t_i^{(e_i)}! \modp & \mbox{if} \; t_i-t_i^{(0)} \; \mbox{is odd} \\  2 (-1)^{i t_i} t_i^{(0)} ! t_i^{(1)} ! \ldots t_i^{(e_i)}! \modp & \mbox{if} \; t_i-t_i^{(0)} \; \mbox{is even} \end{array} \right. .$$

\medskip
Now fix $1 \leq i \leq k$, and let $\{1, \, \ldots , \, e_i \}=I_1^{(i)} \cup I_2^{(i)}$, where $$I_1^{(i)} = \{ j \in \{ 1, \, \ldots , \, e_i \} \, | \, t_i^{(j)}=2k_i^{(j)}+1 \, (k_i^{(j)} \in \N \cup \{ 0 \} ) \}$$and$$I_2^{(i)} = \{ j \in \{ 1, \, \ldots , \, e_i \} \, | \, t_i^{(j)}=2k_i^{(j)} \, (k_i^{(j)} \in \N) \}.$$
We obtain$$\prod_{j=1}^{e_i} \psi_i^{(j)}(1)= \prod_{j \in I_1^{(i)}} 2^{k_i^{(j)}}\chi_{\la_i^{(j)}}(1) \prod_{j \in I_2^{(i)}} 2^{k_i^{(j)}}\chi_{\la_i^{(j)}}(1)= 2^{\sum_{j=1}^{e_i}k_i^{(j)}} \prod_{j=1}^{e_i}\chi_{\la_i^{(j)}}(1).$$
Note that $S_i=|I_1^{(i)}|$ if $\psi_i^{(0)}$ is s.a., while $S_i=|I_1^{(i)}|+1$ if $\psi_i^{(0)}$ is n.s.a.. We thus have
$$2^{\floor{S_i/2}} \prod_{j=1}^{e_i} \psi_i^{(j)}(1) = 2^{\sum_{j=1}^{e_i}k_i^{(j)} + \floor{S_i/2}} \prod_{j=1}^{e_i}\chi_{\la_i^{(j)}}(1) = 2^{\floor{ \sum_{j=1}^{e_i}k_i^{(j)} + S_i/2}} \prod_{j=1}^{e_i}\chi_{\la_i^{(j)}}(1),$$
and
$$\displaystyle \begin{array}{rl}

\displaystyle \floor{ \sum_{j=1}^{e_i}k_i^{(j)} +\frac{ S_i}{2}} & = \displaystyle \floor{ \sum_{j \in I_1^{(i)}} \frac{t_i^{(j)}-1}{2} + \sum_{j \in I_2^{(i)}} \frac{ t_i^{(j)}}{2} + \frac{ S_i}{2} } \\ & =  \left\{ \begin{array}{ll} \displaystyle  \floor{ \sum_{j \in I_1^{(i)}} \frac{t_i^{(j)}}{2} + \sum_{j \in I_2^{(i)}} \frac{t_i^{(j)}}{2} } = \floor{ \frac{t_i-t_i^{(0)}}{2} } & \mbox{if} \; \sa(\psi_i^{(0)})=1 \\ \displaystyle \floor{ \sum_{j \in I_1^{(i)}} \frac{t_i^{(j)}}{2} + \sum_{j \in I_2^{(i)}} \frac{t_i^{(j)} }{2} + \frac{1}{2} } = \floor{ \frac{t_i-t_i^{(0)}+1}{2} } & \mbox{if} \; \sa(\psi_i^{(0)})=-1 \end{array} \right. .  \end{array}$$

We can now compute $D_1^{(+1)}$, $D_0^{(+1)}$, $ D_1^{(-1)}$ and $ D_0^{(-1)}$. Take any $1 \leq i \leq k$, and first suppose that $t_i - t_i^{(0)}$ is odd and $\sa(\chi_{\la_i^{(0)}})=1$. Then $\psi_i^{(0)}$ is n.s.a., $S_i=|I_1^{(i)}| + 1$, $t_i - m_i^{(0)}$ is even and$$\chi_{\la_i^{(0)}}(1)=2^{\floor{ \frac{t_i^{(0)}-m_i^{(0)}}{2} }} \displaystyle \frac{t_i^{(0)}!}{\h(\la_i^{(0)})}.$$
Finally, $ \floor{ \frac{t_i-t_i^{(0)}+1}{2} }= \floor{ \frac{t_i-t_i^{(0)}}{2} } +1$ (since $t_i - t_i^{(0)}$ is odd). We therefore get
$$ \begin{array}{rl} \displaystyle \frac{|T_i^+|}{2^{\floor{S_i/2}-1} \prod_{j=0}^{e_i} \psi_i^{(j)}(1)} & \equiv \pm \frac{ \h(\la_i^{(0)}) \prod_{j=1}^{e_i} h(\la_i^{(j)}) } {2^{ \floor{ \frac{ t_i -t_i^{(0)} +1}{2}} + \floor{ \frac{t_i^{(0)} - m_i^{(0)}}{2}} - 1}} \modp \\ & \equiv \pm \displaystyle \frac{1}{2^{ \floor{ \frac{ t_i -m_i^{(0)} +1 }{2}}-1}} \h(Q_i) \modp \\ & \equiv \pm \displaystyle \frac{2}{2^{ \floor{ \frac{ t_i -m_i^{(0)} }{2}}}} \h(Q_i) \modp
\end{array}$$ 
(since $t_i - m_i^{(0)}$ is even). By similar arguments, we obtain, in all other cases,
$$  \displaystyle \frac{|T_i^+|}{2^{\floor{S_i/2}} \prod_{j=0}^{e_i} \psi_i^{(j)}(1)}  \equiv  \pm \displaystyle \frac{2}{2^{ \floor{ \frac{ t_i -m_i^{(0)} }{2}}}} \h(Q_i) \modp.$$
Finally, we get
$$\frac{|N_G(X)|_{p'}}{\cyc{ (Q_1, \, \ldots , \, Q_k) }(1)_{p'}} \equiv \pm \frac{2^k}{2^{k-1}2^{\floor{s/2}}} \prod_{i=1}^k \frac{1}{2^{\floor{ (t_i - m_i^{(0)})/2 }}} \h(Q_i) \modp ,$$as announced.

\end{proof}

\subsection{Isaacs-Navarro Conjecture}
We can now prove the main result of this section.
\begin{theorem}\label{INC-principal block}
The Isaacs-Navarro Conjecture holds for the principal spin block of $S^+(pw)$ via a sign-preserving bijection.
\end{theorem}
\begin{proof}
Let $B$ be the principal spin block of $G=S^+(pw)$, and $b$ its Brauer correspondent in $N_G(X)$. Let $pw=\sum_{i=1}^k t_i p^i$ be the $p$-adic decomposition of $pw$. By Proposition \ref{height}, $\la \succ pw$ labels a spin character of $B$ of $p'$-degree if and only $\la$ has $\p$-core tower $(R_1^{\la}, \, \ldots , \, R_k^{\la})$ with $|R_i^{\la}|=t_i$ for each $1 \leq i \leq k$. Also, for any such $\la$, we have, by Lemma \ref{tower-sign}, $\sa(\cla)=\sa(\la)=\sa(R_1^{\la}) \ldots \sa(R_k^{\la})$. By the above description of $M(b)$, this implies that 
$$\Phi \colon \left\{ \begin{array}{rcl} M(B) & \longrightarrow & M(b) \\ \cla & \longmapsto & \cyc{ (R_1^{\la}, \, \ldots , \, R_k^{\la})} \end{array} \right.$$
is a sign-preserving bijection. Furthermore, it is immediate from Proposition \ref{degspin} and Proposition \ref{degspinb} that, for any $\cla \in M(B)$,
$$\frac{|G|_{p'}}{\cla (1)_{p'}} \equiv \pm \frac{|N_G(X)|_{p'}}{\cyc{ (Q_1, \, \ldots , \, Q_k) }(1)_{p'}} \modp.$$
This proves the result.

\end{proof}

\section{Main Theorem}
We can now finally give our main theorem:
\begin{theorem}\label{INC}
The Isaacs-Navarro Conjecture holds for all covering groups of the symmetric and alternating groups, whenever $p$ is an odd prime.
\end{theorem}

\begin{proof}
First, let $G$ be any central extension of degree 2 of $S(n)$ or $A(n)$, and $B$ be a $p$-block of $G$. If $B$ is an unfaithful block, then the Isaacs-Navarro Conjecture holds for $B$ by the results of Fong (\cite{Fong}) and Nath (\cite{Nath}). If $B$ is a spin-block of $G$ of weight $w>0$, then the Isaacs-Navarro Conjecture holds for $B$ by Theorem \ref{Reduction Theorem} and Theorem \ref{INC-principal block}. If $w=0$, then $B$ contains a unique spin character (of $p$-defect 0), and the result is immediate.

Finally, the case of the exceptional 6-fold covers of $A(6)$ and $A(7)$ can easily be checked using the character tables given in \cite[6. Appendix]{Michler-Olsson}, or with a computer.

\end{proof}

\noindent\textbf{Acknowledgements.}\quad

%Using Theorem~\ref{nbdecAn},
\bibliographystyle{plain}
\bibliography{referencesJB}

\end{document}